\documentclass[10pt]{amsart}
\usepackage{epsfig}
\usepackage{amssymb}
\usepackage{amscd}
\usepackage{amsthm,amsgen,amsmath,epic,psfrag}
\usepackage[all,2cell,dvips]{xy}\UseAllTwocells\SilentMatrices

\newcommand{\I}{{\mathbb I}}
\newcommand{\R}{{\mathbb R}}
\newcommand{\E}{{{\rm Emb}(\I, M)}}
\newcommand{\Z}{\mathbb{Z}}
\newcommand{\knot}{\mathcal{K}}
\newcommand{\coll}{\mathbf{Co}}
\newcommand{\CO}{\mathcal{C}_4}
\newcommand{\wt}{\widetilde}

\newcommand{\la}{{\langle [}}
\newcommand{\ra}{{] \rangle}}
\newcommand{\Pa}{{\bf Pa}}
\newcommand{\foo}{CM}

\newcommand{\colinear}{collinear}

\newcommand{\colinearity}{collinearity}

\theoremstyle{plain}
\newtheorem{theorem}{Theorem}[section]
\newtheorem{proposition}[theorem]{Proposition}
\newtheorem{lemma}[theorem]{Lemma}
\newtheorem{corollary}[theorem]{Corollary}
\newtheorem{conjecture}[theorem]{Conjecture}
\newtheorem{mydiagram}[theorem]{Figure}
\theoremstyle{definition}
\newtheorem{definition}[theorem]{Definition}

\theoremstyle{remark}
\newtheorem*{remark}{Remark}

\newcommand{\refT}[1]{Theorem~\ref{T:#1}}
\newcommand{\refC}[1]{Corollary~\ref{C:#1}}

\newcommand{\refD}[1]{Definition~\ref{D:#1}}

\DeclareMathOperator{\im}{im}

\begin{document}
\title{New perspectives on self-linking}

\author[R. Budney]{Ryan Budney} 
\address{Department of Mathematics\\
Hylan Building\\
University of Rochester\\
Rochester, NY 14627 
}
\email{rybu@math.rochester.edu}
\author[J. Conant]{James Conant}
\address{Department of Mathematics\\
         Cornell University\\
         Ithaca, NY 14853-4201}
\email{jconant@math.cornell.edu}

\author[K. P. Scannell]{Kevin P. Scannell} 
\address{
Department of Mathematics and Computer Science\\
Saint Louis University\\
Saint Louis, Missouri 63103
}
\email{scannell@slu.edu}

\author[D. Sinha]{Dev Sinha}
\address{Mathematics Department\\
University of Oregon\\
Eugene, OR
97403}
\email{dps@darkwing.uoregon.edu}
\keywords{finite-type invariants, Goodwillie calculus, quadrisecants}
\subjclass{57M27, 55R80, 57R40, 57M25, 55P99}
\thanks{The second author is supported by NSF VIGRE grant DMS-9983660. 
The third author is supported by NSF grant DMS-0072515.}

\begin{abstract}
We initiate the study of classical knots through the homotopy class 
of the $n$th evaluation map of the knot, which is  the 
induced map on the compactified $n$-point configuration space.   
Sending a knot to 
its $n$th evaluation map realizes the space of knots as 
a subspace of what we call the $n$th mapping space model 
for knots.    We compute the homotopy types of the first three
mapping space models, showing that the third model gives 
rise to an  integer-valued
invariant.  We realize this invariant in two ways, in terms of
collinearities of three or four points on the knot, and give some
explicit computations.  We show this invariant
coincides with the second coefficient of the Conway polynomial, thus
giving a new geometric definition of the simplest finite-type invariant.
Finally, using this geometric definition,
we give some new applications of this invariant relating to
quadrisecants in the knot and to complexity of polygonal and 
polynomial realizations of a knot.
\end{abstract}

\maketitle

{\tableofcontents}

\section{Introduction}

Finite-type invariants, introduced by Goussarov and Vassiliev
\cite{Gou94, Va92}, enjoyed an explosion of interest (e.g. \cite{Bi93,BL93,BN95}) 
when Kontsevich, building on work of Drinfeld, showed that they are rationally 
classified by an algebra of trivalent graphs
\cite{Ko93}. The topological meaning of the classifying map 
(defined by the Kontsevich integral) is poorly understood however, and various attempts
have been made to understand finite-type invariants in terms of classical
topology. Bott and Taubes, for example, began the study of finite-type invariants
via the de Rham theory of configuration spaces \cite{BT94, AF97, KT99}. In this
paper we begin to construct universal finite-type invariants \emph{over the
integers}, through homotopy theory and differential topology of 
configuration spaces. As an outcome we
find novel topological interpretations of the first nontrivial finite-type
invariant.

Finite-type invariants in general, and in particular the simplest such invariant 
$c_2$ which is the 
second coefficient  of the Conway polynomial, have been known 
as self-linking invariants for various reasons.
One explanation is that if  $K_+$ and $K_-$ differ by a single crossing 
change, then  $c_2(K_+)-c_2(K_-)  = lk(L_0)$, where $L_0$ is a link formed by 
``resolving" the crossing.  Bott and
Taubes also call $c_2$ a self-linking invariant because their de Rham formula
generalizes the Gauss integral formula for the linking number, as does the
Kontsevich integral.  

\begin{figure}[ht]
$$\includegraphics[width=5cm]{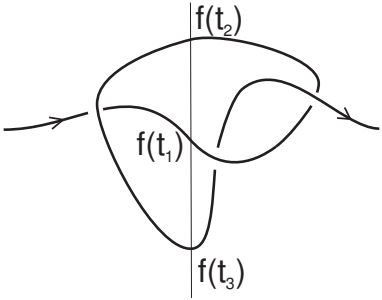}$$
\caption{Collinear points on a knot, giving rise to a point in $\coll_1(\knot)$.} 
\label{co1example}
\end{figure}

We provide a new perspective on self-linking in the following sense.
For technical simplicity, our knots are proper embeddings $f$ of the  interval $\I$
in $\I^3$ with fixed endpoints and tangent vectors at those endpoints,
the space of all such we call ${\rm Emb}(\I, \I^3)$.   
Consider the submanifold $\coll_i(f)$ of ${\rm Int}(\Delta^3)$ consisting of all
$t_1<t_2< t_3$ such that $f(t_1)$, $f(t_2 )$ and $f(t_3)$ are \colinear\  and
such that $f(t_i)$ is between the other two points along the line.  See for
example Figure~\ref{co1example}.

By putting $\knot$ in general position, $\coll_1(f)$ and $\coll_3(f)$ are 
codimension two  submanifolds of ${\rm Int}(\Delta^3)$. By using the appropriate
compactification technology (as developed in Subsection~\ref{compactify} and \cite{Si02a})
they are interiors of 1-manifolds  with  boundary, which we denote $\coll_i[\knot]$,
$i=1,3$, inside the compactification $C_3[\I, \partial]$ of ${\rm Int}(\Delta^3)$.  
Moreover, $\coll_1[\knot]$ and $\coll_3[\knot]$ have boundaries on disjoint  faces of
$C_3[\I, \partial]$, so they have a well-defined linking number. In
Section~\ref{definition}, we define
$\nu_2\colon\pi_0({\rm Emb}(\I,\I^3))\to \Z$ as this linking number of collinearity
submanifolds.
In sections  \ref{examples} and \ref{coquad}, we
show  $\nu_2 = c_2$ thus giving a new geometric interpretation of 
this simplest quantum invariant.

We were led to this self-linking construction as part of a more general study.
In Section~\ref{spaces},  we construct 
an approximating model for the space of embeddings of 
an interval in a manifold, ${\rm Emb}(\I,M)$,
$$C_n: {\rm Emb}(\I,M)\to AM_n(M),$$ 
introduced by the fourth author \cite{Si00} building on the calculus of
embeddings of Goodwillie,  Klein, and Weiss \cite{GW99a,  GKW01, GK00, 
We96, We99}.
When the dimension of $M$ is greater than three, the map $C_n$ induces
isomorphisms on homology and homotopy groups up to degree $n(\dim(M)-3)$. 
For three-manifolds we conjecture a strong relation to finite-type invariants. 
In particular, for $M = \I^3$ we conjecture the following.

\begin{conjecture}\label{conj} The map
$\pi_0(C_n)\colon\pi_0({\rm Emb}(\I,\I^3))\to \pi_0(AM_n(\I^3))$ is a universal
additive type $n-1$ invariant over $\mathbb Z$.
\end{conjecture} 

This conjecture requires explanation.  The connected 
components of ${\rm Emb}(\I,\I^3)$ are knot types, so $\pi_0(C_n)$ is 
indeed a knot invariant. We conjecture that $AM_n(\I^3)$ is always a $2$-fold loop 
space, implying that $\pi_0(AM_n(\I^3))$ is an abelian group. 
Moreover we conjecture that $\pi_0(C_n)$ is
a homomorphism from the monoid of knots, under connected sum, to the
abelian group $\pi_0(AM_n(\I^3))$. Finally, by a \emph{universal type $n-1$
invariant} we mean a knot invariant taking values surjectively
in an abelian group   
such that  every type $n-1$ invariant factors through this map.  

We establish Conjecture~\ref{conj}  for $n \leq 3$.  In section~\ref{components}
we show that 
$AM_1(\I^3)$  is homotopy equivalent to $\Omega S^2$, which is connected, 
and that $AM_2(\I^3)$ is contractible.
We then cite the well-known fact that there are no non-trivial
degree zero or one knot invariants.  We focus on $n=3$, where our work is grounded in 
computations similar to those of \cite{SS02}, in which the last
two authors compute rational homotopy groups of spaces of knots in  
even-dimensional Euclidean spaces.  For classical knots, they
found that $\pi_0(AM_3(\I^3))$ is isomorphic to $\pi_3$ of the homotopy
fiber $F$ of  the inclusion $S^2\vee S^2\to S^2\times S^2$.  In fact, in 
section~\ref{components} we establish that $AM_3(\I^3) \simeq \Omega^3 F$.  
Classically, $\pi_3(F)$ is well known to be isomorphic to the integers, with 
the isomorphism given  by a linking number or Hopf invariant $\mu_2$. We show  
that $\mu_2 \circ \pi_0(C_n) = \nu_2 = c_2$, the first non-trivial Vassiliev invariant. 
Higher $n$ will be considered in \cite{BCSS03}, in which we plan to define the
knot invariants given by $\pi_0(C_n)$ through linking invariants of what we
call ``coincidence submanifolds'' in the parameter space $C_n[I, \partial]$. 
Constructing knot invariants through Hopf invariants of
natural submanifolds of this parameter space is our new perspective on self-linking.

Further evidence for Conjecture~\ref{conj} is the thesis of I.~Volic,
which shows that the 
Bott-Taubes invariants factor through $\pi_0(C_{2n})$ \cite{Vo03}. 
Since Bott-Taubes invariants generate the rational vector space of finite 
type invariants, this shows that $\pi_0(C_{2n})$ rationally classifies finite 
type invariants. Our approach is complementary in that it is over the integers, 
and moreover leads to novel geometric
consequences, even at the lowest non-trivial degree.  Another approach 
to finite-type knot theory which
may be closely related to ours is that of arrow diagrams \cite{GPV00}, which 
following ideas of  \cite{PV98} might be understood through composing
the evaluation map of this paper with projections to products of spheres.
We hope to explore this connection in future work.

For general three-manifolds $M$ 
the theory of finite-type invariants is not well understood, and accordingly
neither is $AM_n(M)$, whose components do not seem to have an additive
structure.  Nonetheless, we suspect that Conjecture~\ref{conj} will have an
analogue in this setting.
It may be helpful to restate the conjecture dually. Define two knots to be 
$n-1$-equivalent if  they share all type $n-1$  invariants. 
As proven by Goussarov \cite{Gou94}, knots modulo
$n-1$ equivalence form an abelian group under connected sum. One can 
formulate $n-1$ equivalence in other ways, as equivalence
up to:
\begin{itemize}
\item Tying a pure braid in the $n$th term of the lower central series of the 
pure braid group into some strands of a knot. \cite{St98}
\item Simple clasper surgeries of degree $n$. \cite{Ha00}
\item Capped grope cobordism of class $n$. \cite{CT03}
\end{itemize}
Conjecture~\ref{conj} can be restated as
saying $\pi_0(C_n)$ induces an isomorphism
$$\left(\pi_0({\rm Emb}(\I,\I^3))/ {\rm (n-1)-equivalence}\right) 
        \cong \pi_0(AM_n(\I^3)),$$
and we suspect that  different points of view on $n-1$ equivalence
will be helpful in  studying analogues of Conjecture~\ref{conj} for 
arbitrary three-manifolds.

One of the surprises we encountered during this investigation was
that our invariant has a natural interpretation in terms of counting
quadrisecants, which by definition are collinearities of four points along the
knot.  Quadrisecants have appeared previously in knot theory 
\cite{Pa33, MM82, Ku94}. In
Section~\ref{coquad}, we explain this connection in detail.
This quadrisecant interpretation leads to lower bounds on
stick number for polygonal knots and the degree of polynomial 
needed to represent a knot 
in terms of $c_2$,  also explained in  section~\ref{coquad}. 

{\bf Acknowledgments} We wish to thank Tom Goodwillie, 
Louis Kauffman, John Sullivan, and Arkady Vaintrob and the referee 
for their interest and suggestions.

\section{Spaces of knots and evaluation maps}\label{spaces}

The appropriate construction of an evaluation map for knots is our connection 
between algebraic and geometric topology.  
Evaluation maps have long been a tool in the
study of links.  The linking number of a two component link can be 
understood as the homotopy class of the evaluation map from the space 
of two points on the link, one on each component, to $C_2(\R^3)$, 
the space of ordered pairs in
$\R^3$. Up to homotopy, this is a map from $S^1\times S^1$ to
$S^2$, which is thus characterized by degree.  
This degree is computed  by the Gauss integral.  In fact, Koschorke has 
shown \cite{Ko97} that all Milnor linking
numbers may be understood through evaluation maps from $(S^1)^k$ to
$C_k(\R^3)$.  For knots, the situation is more subtle because one is led to
consider maps from the space of configurations on the knot, which is an open
simplex, to $C_k(\R^3)$.  
The technical heart of the matter is what to do ``on the boundary'' of
the configuration space, including what boundary to use in the first place.  In
\cite{Si00} the fourth author showed that the appropriate boundary conditions
are prescribed when one relates the calculus of embeddings to the evaluation map.  

\subsection{Compactifications of configuration spaces}\label{compactify}

Key technical tools are the Fulton-MacPherson compactifications of configuration spaces.  
Readers familiar with these constructions may skip 
most of this section, but should familiarize themselves with our labeling 
scheme for strata.  For manifolds, there are two versions of these compactifications, 
a projective version which can be defined by 
an immediate translation  of the Fulton-MacPherson
construction, and a closely related spherical version which was first
constructed by Axelrod and Singer \cite{AS94}.  
We use the spherical version, in many ways most natural for manifolds,
which up to homeomorphism is just the space $M^n \setminus N$, 
where $N$ is a ``tubular''
neighborhood of the fat diagonal.  These compactifications have many
properties which are not immediate from the $M^n \setminus N$ model.  See
\cite{Si02a} for a full development of these compactifications using the
following simple construction similar to one given for Euclidean spaces
by Kontsevich \cite{Ko99} and Gaiffi \cite{Ga03}.  

\begin{definition}
Let $\pi_{ij} \colon C_n(\R^k) \to S^{k-1}$ be the map which sends 
${\bf x} = (x_1, \ldots, x_n)$ to the unit vector in the direction of $x_i - x_j$.
Let $s_{ijk} \colon C_n(\R^k) \to (0, \infty) \subset \I$, where the latter
inclusion is defined through the arctangent function, 
be the map which sends ${\bf x}$ to
$\left(|x_i - x_j|/|x_i - x_k|\right)$.  Let $B_{n,k}$ be the product 
$(\R^k)^n \times (S^{k-1})^{\binom{n}{2}} \times \I^{\binom{n}{3}}$.
\end{definition}

\begin{definition} 
Define $C_n[\R^k]$ to be the closure of the image of $C_n(\R^k)$ under the
 map $\iota \times \prod \pi_{ij} \times \prod s_{ijk}$ to $B_{n,k},$ where 
$\iota$ is the standard inclusion of $C_n(\R^k)$ in $(\R^k)^n$.  Define 
$C_n[M]$ for a general  manifold
$M$ by embedding $M$ in some $\R^k$ in order to define the restrictions of the 
maps $\pi_{ij}$ and $s_{ijk}$, and then taking the closure of $C_n(M)$ in 
$M^n \times \prod (S^{k-1}) \times \prod \I$.  For $M = \I^k$ we use the standard 
embedding of $\I^k$ in $\R^k$.
\end{definition}

These closures have the following important properties.

\begin{theorem}\label{T:compact} (see \cite{Si02a})
\begin{itemize}
\item $C_n[M]$ is a manifold with corners whose interior is $C_n(M)$, and 
which thus has the same homotopy type as $C_n(M)$.  It is independent
of the embedding of $M$ in $\R^k$, and it is compact when $M$ is.
\item The inclusion of $C_n(M)$ in $M^n$ extends to a surjective map $p$ from
$C_n[M]$ to $M^n$ which is a homeomorphism over points in $C_n(M)$.
\item If, in the projection of $x \in C_n[M]$ onto $M^n$, some $x_i =
x_j$, then there is a well-defined $v_{ij} \in STM$ sitting over $x_i$, 
which we call the {\em relative vector} and which gives the 
``infinitesimal relative position'' of $x_i$ and $x_j$.  
For $M = \R^k$, this vector is given by the extension of the map $\pi_{ij}$, which 
is in turn the restriction of the projection of $B_{n,k}$ onto the $ij$ factor of 
$S^{k-1}$.
\item If $f \colon M \to N$ is an embedding the induced map on open
configuration spaces extends uniquely to a
map, which we call the evaluation map,
$C_n[f] \colon C_n[M] \to C_n[N]$ which preserves the stratification of these
spaces as manifolds with corners.  This map is defined by choosing
the ambient embedding of $M$ in $\R^k$ to be the composite of $f$ with 
the ambient embedding of $N$.
\end{itemize}
\end{theorem}

The best understanding of these compactifications comes from their
stratification.  When $M$ is Euclidean space, these
strata are simple to describe once we develop the appropriate combinatorics to
enumerate them. Here, we choose to develop this combinatorics in terms of
parenthesizations.

\begin{definition}
A {\em (partial)  parenthesization} $\mathcal{P}$ of a set $T$ is an unordered
collection $\{ A_i \}$ of subsets of $T$, each having cardinality at least two,
such that for any $i,j$ either $A_i \subset A_j$ or $A_i$ and $A_j$ are disjoint.  
We denote a parenthesization by a nested listing of elements of the $A_i$ using 
parentheses and equal signs.
\end{definition}

For example, $(3=4), ((1=2)=6)$ represents a parenthesization of $\{1, \ldots, 6\}$ 
whose subsets are $\{3,4\}$, $\{1,2\}$ and $\{1,2,6\}$.

\begin{definition} 
Let $\Pa(T)$ be the set of parenthesizations of $T$.  Define an ordering on $\Pa(T)$ 
by $\mathcal{P} \leq \mathcal{P}'$ if $\mathcal{P} \subseteq \mathcal{P}'$.
A {\em total parenthesization} is a maximal element under this ordering.
\end{definition}

For example, $\Pa(\{1,2,3\})$ is given as follows.

\[
\xymatrix {
{} & ((1=2)=3) & ((2=3)=1) & ((3=1)=2) & {} \\
{} & (1=2=3) \ar[u] \ar[ur] \ar[urr] & (1=2) \ar[ul] & (2=3) \ar[ul] &
(3=1) \ar[ul] \\ {} & {} & \emptyset \ar[ul] \ar[u] \ar[ur] \ar[urr] & {}
& {} }
\]


In our applications we need the combinatorics of maximal subsets within
a parenthesization.

\begin{definition}  
From a parenthesization $\mathcal{P}$ of $T$ we form a sequence
$\lambda_i(\mathcal{P})$ which for $i > 0$ are subsets of $\mathcal{P}$ 
by setting
$\lambda_0(\mathcal{P}) = \{ T \}$ and inductively the elements of
$\lambda_{i+1}(\mathcal{P})$ are maximal elements in the ordering given by
inclusion of subsets of $T$ (that is, largest subsets), 
of $\mathcal{P} \setminus \bigcup_{k=0}^i \lambda_k(\mathcal{P})$.

Define an equivalence relation $\sim_i$ on each $Q \in \lambda_i(\mathcal{P})$ 
by $x \sim_i y$ if $x, y \in Q'$ for some $Q' \in \lambda_{i+1}(\mathcal{P})$.  
For such a $Q$, let $|Q|$ be the number of equivalence classes of $Q$ under 
$\sim_i$.
\end{definition}

Recall that the strata of a manifold with corners form a poset by the relation $S <
T$ if $S$ is in the closure of $T$.  We are now ready to describe the strata of
$C_n[\R^k]$.

\begin{definition}
Let $\wt{C}_n(\R^k)$ be the quotient of $C_n(\R^k)$ by the action of scaling all
points by positive real numbers and by  translation of all points.
\end{definition}

So for example $\wt{C}_2(\R^k)$ is diffeomorphic to $S^{k-1}$.

\begin{theorem}[See \cite{Si02a}]\label{T:strata}
If $M$ has no boundary, the poset of strata of $C_n[M]$ is isomorphic to $\Pa([n])$,  
where $[n] = \{1, \ldots, n\}$.  Given a parenthesization $\mathcal{P}$ of $[n]$, the
codimension of the corresponding stratum $S_{\mathcal{P}}$ is the cardinality of
$\mathcal{P}$.  When $M = \R^k$, the stratum $S_{\mathcal{P}}$  is diffeomorphic to  
$$ C_{|T|}(\R^k) \times \prod_{i > 0} \prod_{Q \in \lambda_i(\mathcal{P})} 
                     \wt{C}_{|Q|}(\R^k).$$
\end{theorem}

For $M=\R^k$, we give an informal indication of how the codimension zero 
stratum  $C_n(\R^k)$ and the codimension one strata 
$C_{n-\ell+1}(\R^k) \times \wt{C}_\ell(\R^k)$ are topologized together as follows.  
Let $\{\bf{x}_i = (x_{i1}, \ldots, x_{in})\}$ 
be a sequence of points in $C_n(\R^k)$, and consider its image
in $(\R^k)^n$.  For $\{\bf{x}_i\}$ to converge to a point in   
$C_{n-\ell+i}(\R^k) \times \wt{C}_\ell(\R^k)$, its image in $(\R^k)^n$ must
converge to a point ${\bf x}_\infty = (x_{\infty 1}, \ldots, x_{\infty n})$ 
which has $\ell$ coordinates which coincide, say 
$x_{\infty j_1} = \cdots = x_{\infty j_\ell}$.   The other
coordinates, $x_{\infty i}$ for $i \notin \{j_k\}$, along with say 
$x_{\infty j_1}$ define a point in $C_{n-\ell+1}$.  
Moreover, to obtain a point in $\wt{C}_\ell(\R^k)$ in the limit from 
$\{\bf{x}_i\}$ we must be able to take the limit of 
$\{ (x_{i j_1}, \ldots, x_{i j_\ell}) / \sim \}$ in $\wt{C}_\ell(\R^k)$.

\medskip

In our variant of knot theory, the ambient manifold $\I^3$
comes equipped with two points on its boundary to which the endpoints of the
interval must map.  

\begin{definition}
If $M$ is a manifold with boundary with two distinguished points $y_0$
and $y_1$ in its boundary,  define $C_n[M, \partial]$ to be the subspace of
$C_{n+2}[M]$ of points whose images under $p$ in  $M^{n+2}$ are of the form
$(y_0, x_1, \ldots, x_n, y_1)$. 
\end{definition}

One case of interest is of course when $M$ is an interval and the two
distinguished points are its endpoints.  The configuration space of points in
the open interval is a union of open simplices.  At times, we will compactify 
the open simplex by embedding it in the standard way in the closed simplex.  
We call $\Delta^n=\{ {\bf t} = 
(t_1,\ldots,t_n) : 0 = t_0  \leq t_1 \leq \cdots \leq t_n \leq t_{n+1} = 1\}$  
the {\em naive compactification} of configurations in the interval and
note that our projection map $p$ sends $C_n[\I,
\partial]$ to $\Delta^n$.  The compactification $C_n[\I, \partial]$, which has
arisen elsewhere in topology \cite{St70}, 
is a polyhedron with a beautiful combinatorial
description.

\begin{definition}
A subset $A$ of a totally ordered set $S$ is {\em consecutive} if $x,z \in A$
and $x < y < z$ implies $y \in A$, and $A$ has two or
more elements.  An  {\em ordered parenthesization} of $S$ is
one in which all of the subsets in the parenthesization are consecutive and
proper.
\end{definition}

\begin{definition}
The $n$th {\em associahedron} (or {\em Stasheff polytope}), 
$A_n$, is the polytope whose
faces $A_n^P$ are labeled by the ordered parenthesizations of
$\{0, \ldots, n+1\}$ and where $A_n^P$ is a face of $A_n^{P'}$ if $P \subset P'$.
\end{definition}

Alternatively, the barycentric subdivision of $A_n$ is isomorphic to
the nerve (also known as the realization or the order complex) of the 
poset of ordered parenthesizations of $\{0, \ldots, n+1\}$.  The name
associahedron comes from the fact that the vertices are labeled by total ordered
parenthesizations of $\{0, \ldots, n+1\}$, 
which may be viewed as ways to associate a word with
$n$ letters.

\begin{proposition}[See \cite{Si02a}]\label{P:stash}
Each component of $C_n[\I, \partial]$ is isomorphic
as a manifold with corners to the $n$th associahedron $A_n$.  
\end{proposition}


\begin{center}
\begin{minipage}{7cm}\begin{mydiagram}\label{a3pic}
\begin{center}
{\psfrag{x}{$t_1$}\psfrag{y}{$t_2$}\psfrag{z}{$t_3$}
$$ \includegraphics[width=8cm]{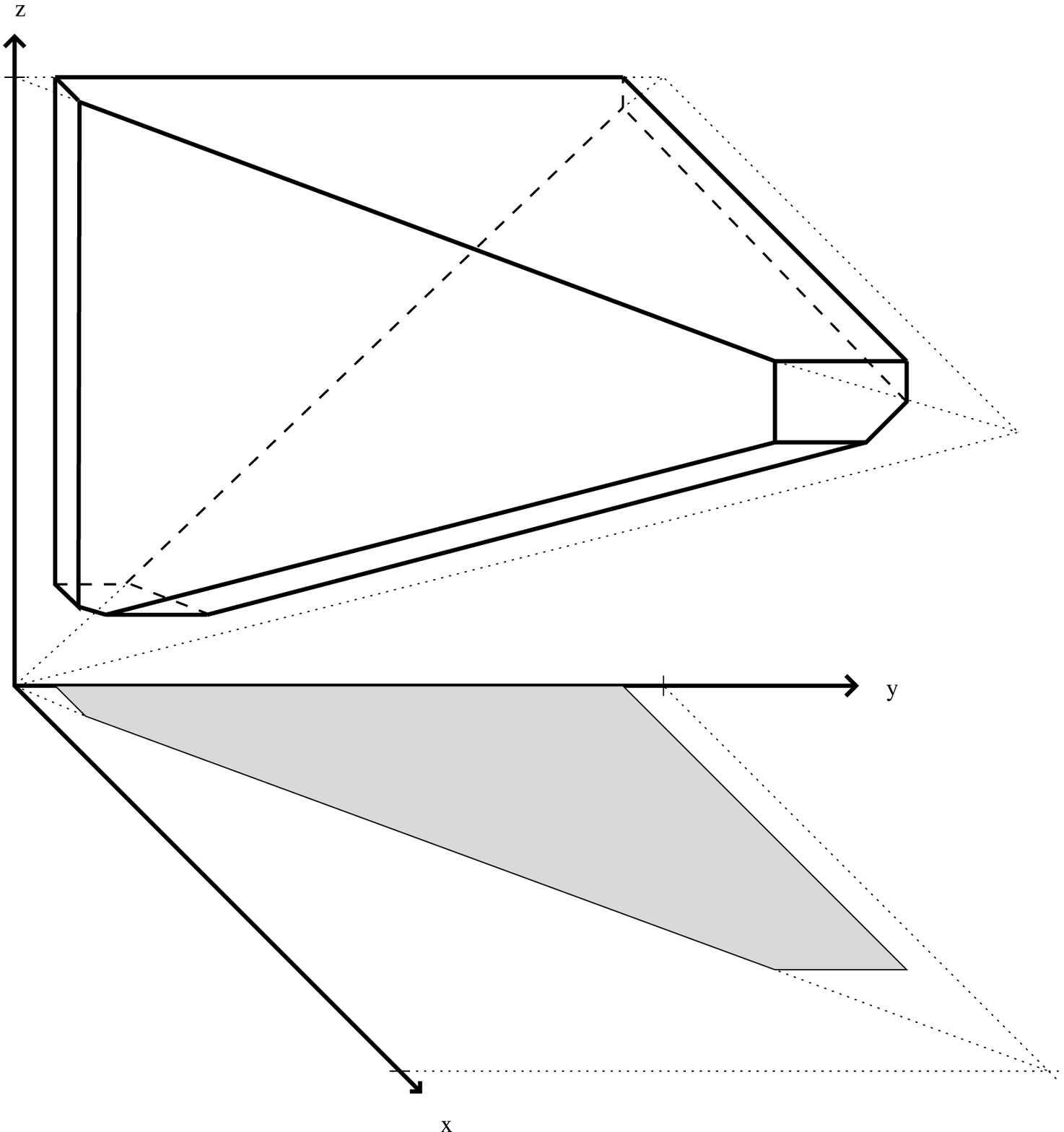}$$}
The $3$rd Stasheff polytope $A_3 = C_3[\I,\partial]$
\end{center}\end{mydiagram}\end{minipage}
\end{center}

Finally, we include tangent vectors in the picture as follows.
\begin{definition}
Let $C'_n[M]$ be defined by the pull-back square
$$
\begin{CD}
C'_n [M]  @>>> (STM)^n \\
@VVV   @VVV \\
C_n[M]  @>p>> M^n.
\end{CD}
$$
If $M$ has two distinguished points on its boundary with one inward-pointing 
and one outward-pointing tangent  vector at those points, 
let $C'_n[M, \partial]$ be the subspace of 
$C'_{n+2}[M]$ of  points whose first and last projection onto $STM$ are 
given by those  distinguished tangent vectors.
\end{definition}

\subsection{The mapping space model}\label{mapping}

We are ready to define our mapping space model for knot spaces.  
As we have mentioned (Theorem \ref{T:compact}), the evaluation map of a knot
$f$, namely $C_n[f] \colon A_n = C_n[\I,\partial] \to C'_n[\I^3, \partial]$,
is stratum preserving.  But note as well that for a point ${\bf t}$  
on the boundary of $A_n$ such that $p({\bf{t}})$ has $t_i = t_{i+1} = t_{i+2}$, 
$C_n[f]({\bf t})$ will, in the decomposition of \refT{strata},
project onto a point in $\wt{C}_3(\R^k)$ which is
degenerate in the sense that all points are aligned along a single vector,
namely the tangent vector to the knot at $f(t_i)$.

\begin{definition}
A point $x \in C_n[M]$ is {\em aligned} if (naming its image in $M^n$ by 
$(x_1, \ldots, x_n)$) for each collection $x_i = x_j = x_k = \cdots$, the
relative vectors $v_{ij}$ with $i<j$ are all equal. A point in $C'_n[M]$
is {\em aligned} if its projection onto $C_n[M]$ is 
aligned and, if we say $x_i = x_j$ as above, the relative vector
$v_{ij}$ serves as the tangent vector at both $x_i$ and $x_j$.
\end{definition}

The subspace of a stratum consisting of points which are aligned
is called the {\em aligned sub-stratum} of that stratum.

\begin{definition}
A stratum-preserving map from $A_n$ to $C'_n[M, \partial]$ is {\em aligned}
if points in its image are aligned.  Let $AM_n(M)$ denote the space of
aligned stratum-preserving maps from $A_n$ to $C'_n[M, \partial]$.
Let $C_n \colon \E \to AM_n(M)$ be the map which sends 
a knot $f$ to its evaluation map $C_n[f]$.
\end{definition}

In \cite{Si00} the fourth author establishes the following.

\begin{theorem}\label{T:modeleq}
The map $C_n \colon \E \to AM_n(M)$ is a model for the $n$th degree
approximation to $\E$ in the calculus of embeddings.
\end{theorem}

As mentioned in the introduction, this result implies that the evaluation
map captures all the topology of families of knots in ambient spaces of
dimension greater than three.  Thus, we are motivated to understand the 
ramifications of this map in knot theory,
and in particular are interested in the effect of $C_n$ on components of 
${\rm Emb}(\I, \I^3)$.  

\section{The homotopy types of the first three mapping space
models}\label{components}

\nocite{FH01}

In this section, we examine the homotopy types of the first three mapping 
space models $AM_i(\I^3)$ for $i = 1,2,3$.  The computations of these
homotopy types can be determined quickly if one cites \refT{modeleq}
and builds on the analysis
of the embedding tower for knots of Goodwillie and Weiss in section five
of \cite{GW99a}, or alternately if one uses the cosimplicial
model for this embedding tower as first envisioned by Goodwillie and developed
in \cite{Si00}.   Indeed, the main computations of this section were previously
known to Goodwillie and Weiss.  We choose a more self-contained 
approach so that the method of analysis fits with our geometric development
of the self-linking invariant  in section~\ref{definition}.

One may view $AM_n(\I^3)$ as a 
multi-relative mapping space.  Standard fibration sequences in algebraic topology 
allow one to build such a mapping space with building blocks which are loop 
spaces.  In this section we employ some standard facts about categories of 
diagrams of spaces in order to organize and simplify such an analysis of
$AM_n(\I^3)$.

First note that $A_1 = \I$, $C'_1[\I^3] \simeq S^2$ and an aligned map from 
$A_1$ has fixed boundary points, so that $AM_1[\I^3] \simeq \Omega S^2$,
the space of based maps from a circle to $S^2$.  Moreover, the map from
the space of knots sends $f$ to $u(f)$ where $u(f) (t) = \frac{f'(t)}{||f'(t)||}$.
This identification is consistent with the fact that in the calculus of embeddings
the first approximation to a space of embeddings is the corresponding
space of immersions, which in turn can be identified through the 
Hirsch-Smale theorem \cite{We96}.  Because $AM_1[\I^3]$ is connected, 
$C_1$ does not give rise to any knot invariants.

We proceed to show that the second mapping space model $AM_2(\I^3)$ is 
contractible, from which we conclude that $C_2$ does not give rise to
any knot invariants either.  This is an analogue of the fact that there are no
non-trivial degree zero or one finite-type invariants.
Finally, we show that $AM_3(\I^3)$ is homotopy equivalent to the third loop 
space of  the homotopy fiber of the inclusion 
$S^2 \vee S^2 \hookrightarrow S^2 \times S^2$.  
By enumerating components, we show that the evaluation map
$C_3$ can give rise to an integer-valued knot invariant, which is
the invariant of study in this paper, defined in Section~\ref{definition}.  

For the purposes of the differential topology, it is essential to use the Fulton-
MacPherson compactifications of configuration spaces.  For the purposes of the 
algebraic topology in this section, however,  it is helpful
to use a variant of these compactifications called the 
{\em simplicial variant}, essentially first defined by Kontsevich \cite{Ko99}
and fully developed in \cite{Si02a}. 
One main advantage is that, when dealing with the configuration space
of $n$ ordered points in the interval, this compactification coincides with the
naive compactification $\Delta^n$.

\begin{definition}
For a manifold $M$ with embedding in $\R^k$ let 
$C_n \la M \ra$ be the closure of the image of $C_n(M)$ in
$A_{k,n}(M) = (M)^n \times (S^{k-1})^{\binom{n}{2}}$ under the map $\iota
\times \prod \pi_{ij}|_M$.  Let $\pi_{ij}$ denote the extension of the maps of 
the same name on the interior of $C_n \la M \ra$ which are defined simply by 
projecting $A_{k,n}(M)$ onto the $ij$ factor of $S^{k-1}$.  
For a manifold with two distinguished  points on its boundary, 
define $C_n \la M, \partial \ra$ and $C'_n \la M, \partial \ra$
as subspaces of $C_{n+2} \la M \ra$ as before.
\end{definition}

For these compactifications, the combinatorics of the strata which arise
in the image of evaluation maps are simplicial.  We formalize this stratification 
as follows.  


\begin{definition}
Define a {\em \mbox{$\Delta_n$-space}} to be a 
collection of spaces $X_S$, one for each
$S \subsetneq \{0, \ldots, n\}$, with a map $f_ {S' \subseteq S} : X_{S} \to X_{S'}$
for each pair $S' \subseteq S$ such that $f_{S'' \subseteq S'} \circ f_{S' \subseteq S}
= f_{S'' \subseteq S}$.
\end{definition}

In other words, a $\Delta_n$ space is a contravariant functor from the 
category of proper subsets of $\{0, \ldots, n\}$ to the category of spaces.
In our examples, all $f_{S \subseteq S'}$ will be inclusions (cofibrations in fact -
see below), so that the
$\Delta_n$-structure describes a stratification.
We often use the notation $X_\bullet$ for \mbox{$\Delta_n$-spaces}. 
As one might expect, there is a natural way of making $\Delta^n$ into a 
\mbox{$\Delta_n$-space}.  Namely, take $\Delta_\emptyset = \Delta^n$
and $\Delta_S$ to be the face for 
which $t_i= t_{i+1}$ if $i \in S$ (recall that $t_0 = 0$ and $t_{n+1}=1$
by convention).  
We write $\Delta^n_\bullet$ for this \mbox{$\Delta_n$-space}.

\begin{definition}
Define a  \mbox{$\Delta_n$-space} $C^n_\bullet$ as follows.
\begin{itemize}
\item $C^n_\emptyset = C_n'\la \I^{3}, \partial \ra = 
          C_n \la \I^{3}, \partial \ra \times (S^2)^n$.  
\item  $C^n_S$ is the subspace of all 
$(x_j) \times (v_k) \in C_n'\la \I^{3}, \partial \ra$ such that for each $i \in S$, 
$x_{i} = x_{i+1}$ and $\pi_{i,i+1}((x_j)) = v_i = v_{i+1}$.
\end{itemize}
\end{definition} 

It is straightforward to check that  $C^n_\bullet$ is a 
\mbox{$\Delta^n$-space} with  $C^n_S$  homeomorphic to 
$C_{n - \#S} '\la \I^3, \partial \ra$.   Moreover if we replace
$\I^3$ by $\I$ in the definition of $C^n_\bullet$, we recover
$\Delta^n_\bullet$.

Many constructions for spaces have immediate analogues for 
\mbox{$\Delta_n$-spaces}.  A map between \mbox{$\Delta_n$-spaces} $X_\bullet$ 
and $Y_\bullet$ is simply a collection of compatible maps from each $X_S$ to $Y_S$.  Thus, the set of all maps 
between two \mbox{$\Delta_n$-spaces} is naturally a subspace
of $\prod_S{\rm Maps}(X_{S}, Y_{S})$.  The case in which 
$X_\bullet = \Delta^n_\bullet$ is of particular use.  We call the space 
of maps from $\Delta^n_\bullet$ to $Y_\bullet$ the 
{\em corealization} of $Y_\bullet$ and denote it $|Y_\bullet|$. 

\begin{remark}
The corealization of $Y_\bullet$ is homeomorphic to the homotopy limit of 
$Y_\bullet$.  The language of homotopy limits is used in \cite{Si00}.
\end{remark}

The following theorem shows that, for the purposes of this section,
the corealization $|C^n_\bullet|$ will replace the mapping space 
model $AM_n(\I^3)$.  We will use the notation $\foo_n$ for 
$|C^n_\bullet|$ and call it the {\em $n$th \mbox{$\Delta$-space} model}.

Consider the evaluation map for a knot using the cosimplicial compactification 
$C_n \colon C_n \la \I, \partial \ra = \Delta^n \to C_n' \la \I^3, \partial \ra$.  
Because this evaluation map is aligned and stratum preserving, it defines a point 
in $\foo_n$. The following can be deduced from Theorem 6.2 in \cite{Si00}.

\begin{theorem}
The mapping space model $AM_n(\I^3)$ is homotopy equivalent to
$\foo_n$. The maps from ${\rm Emb}(\I, \I^3)$ to these mapping
spaces given by evaluation maps agree in the homotopy category.
\end{theorem}

In light of this theorem, we will now  analyze the homotopy types of 
$\foo_2$ and $\foo_3$.  
The notion of a fibration of \mbox{$\Delta_n$-spaces} is helpful;
this is simply a map of \mbox{$\Delta_n$-spaces} with 
the property that each $X_S \to Y_S$ is a (Serre) fibration 
in the usual sense.  For any $S \subset \{0, \ldots, n\}$ let $\overline{S}$ be
its complement.   We say that a $\Delta_n$-space
is cofibrant if each $X_S$ is a CW complex, the structure maps are
inclusions of  subcomplexes, and 
$X_{\overline{S}} \cap X_{\overline{S'}} = X_{\overline{S \cap S'}}$, where by
convention $X_{\{0, \ldots, n\}}$ is empty.  
Finally, if $A$ is a space and $X_\bullet$
is a $\Delta_n$-space let $A \times X_\bullet$ be the $\Delta_n$
space whose $S$ entry is $A \times X_S$ and whose structure
maps are the product of the identity on $A$ with those for $X_\bullet$.

The following lemma is a standard part of developing an enriched model
structure for $\Delta_n$-spaces
(as can be developed for any diagram category) \cite{DK83, Hi03}.

\begin{lemma}\label{L:fibr}
Let $X_\bullet$ be a cofibrant \mbox{$\Delta_n$-space} such that each 
$X_S$ is a finite CW complex. Then a fibration of 
\mbox{$\Delta_n$-spaces} $Y_\bullet \to Z_\bullet$ induces
a fibration on the mapping spaces ${\rm Maps}(X_\bullet, Y_\bullet) \to
{\rm Maps}(X_\bullet, Z_\bullet)$.
In particular, a fibration of \mbox{$\Delta_n$-spaces} gives rise to a fibration
of their corealizations.
\end{lemma}

\begin{proof}
To establish the lifting property
we must show that if $A \hookrightarrow B$ is a cofibration
and a weak equivalence and we are given a square
$$
\begin{CD}
A @>>> { \rm Maps}(X_\bullet, Y_\bullet)\\
@VVV  @VVV \\
B @>>> {\rm Maps}(X_\bullet, Z_\bullet),
\end{CD}
$$
then there is a map $B \to { \rm Maps}(X_\bullet, Y_\bullet)$ which
commutes with the maps in this square.  By adjointness
it suffices to construct a lift of a map $B \times X_\bullet \to Z_\bullet$
to $Y_\bullet$ extending a lift given on $A \times X_\bullet$.
We construct this lift inductively over the $B \times X_S$ starting with
$S$ with $\#S = n$, for which we simply cite that $A \times X_S \hookrightarrow
B \times X_S$ is a cofibration and weak equivalence and use
the fact that $Y_S \to Z_S$ is a fibration.  Suppose we have
constructed a lift on all $B \times X_{S'}$ for $S \subset S'$.  
From the definition of cofibrant $\Delta_n$-space it follows
that the map from the colimit of the $X_{S'}$ over $S \subset S'$ 
to $X_S$ is a cofibration.
Therefore, the union of $A \times X_S$ with 
$\bigcup_{S \subset S'} (B \times X_{S'})$ sits in $B \times X_S$ as a
cofibration and a weak equivalence, so we may again extend
by the fact that $Y_S \to Z_S$ is a fibration.
\end{proof}

Our basic strategy in analyzing $\foo_2$ and $\foo_3$ 
is to fiber $C^n_\bullet$ as 
$F^n_\bullet \to C^n_\bullet \overset{\pi}{\to} K^n_\bullet$ where $K^n_\bullet$
is of a form in which its corealization can be shown to be contractible
explicitly.  For $n=2$, $\pi$ will be an equivalence, and for $n=3$ we 
identify $|F^3_\bullet|$.

\begin{theorem}\label{T:amtwo}
The space $\foo_2$ is contractible.
\end{theorem}

\begin{proof}
Our knots  begin at the ``north'' side of $\I^3$ and end at the ``south'' side, 
iso that $v_0(x) = v_3(x) = \ast$, where $\ast$ denotes 
the south pole of $S^2$.
We make use of the degree one map 
$\sigma : S^2 \to S^2$ which collapses the southern hemisphere to
the south pole, which gives
$\sigma \circ \pi_{0,i} = \sigma \circ \pi_{i,3} = \ast$ for $i=1,2$.
Define $K^2_\bullet$, which is obtained from $C^2_\bullet$ basically by
projecting $C_2 \la \I^3, \partial \ra$ onto $S^2$ through 
$\sigma \circ \pi_{12}$. 
Explicitly, we have the following:
\begin{align*}
K^2_{\emptyset} = (S^2)^3  && K^2_{\{a\}} = S^2 && K^2_{\{a,b\}} = pt.\\
i_{\{a\} \subset \{a,b\}}(pt.) = \ast && i_{\emptyset \subset \{2\}}(x) = (\ast,x,\ast) &&
i_{\emptyset \subset \{0\}}(x) = (\ast,\ast,x) && i_{\emptyset \subset \{1\}}(x) = 
        (\sigma(x),x,x)
\end{align*}

We view a \mbox{$\Delta^n$-space} as a diagram of spaces, 
one sitting at each barycenter in the $n$-simplex, with a morphism for
each face relation.  
As a diagram, $K^2_\bullet$ is the following:

\[
\xymatrix @R=3ex @C=-2ex {
{} & {} & K^2_{\{0,2\}} = \ast \ar[dddl] \ar[dddr] & {} & {} \\
{} & {} & {} & {} & {} \\
{} & {} & {} & {} & {} \\
{} & K^2_{\{0\}} = S^2 \ar[rd] & {} & K^2_{\{2\}} = S^2 \ar[ld] & {} \\
{} & {} & K^2_\emptyset = (S^2)^3 & {} & {} \\
{} & {} & {} & {} & {} \\
K^2_{\{0,1\}} = \ast \ar[ruuu] \ar[rr] & {} & K^2_{\{1\}} = S^2 \ar[uu]
& {} & K^2_{\{1,2\}} = \ast \ar[ll] \ar[luuu] }
\]

The fibration $\pi$ from $C^2_\bullet$ to
$K^2_\bullet$ is defined as $(\sigma \circ \pi_{12}) \times v_1
\times v_2$, a map of \mbox{$\Delta_2$-spaces}.
Furthermore, the fiber of this map restricted to each
$C^2_S$ is contractible, so the induced map from $|C^2_\bullet|$ to 
$|K^2_\bullet|$ is an equivalence.  We now show that $|K^2_\bullet|$ is
contractible.

Define a \mbox{$\Delta_2$-space} $Y_\bullet$ by projecting 
$K^2_{\emptyset}$  onto its first factor of $S^2$, and taking  $Y_S$ to be the 
image of $K^2_S$ under this projection.   
In particular $Y_{\{2\}} = Y_{\{0\}} = \ast$.
The projection from $K^2_\bullet$ to $Y_\bullet$ is readily seen to be a fibration 
of \mbox{$\Delta_2$-spaces}.  It is also easy to see that $|Y_\bullet|$ is 
contractible; unraveling definitions, it is the space of maps from 
$\Delta^2$ to $S^2$ sending  the $\{0\}$ and $\{2\}$ sides of $\Delta^2$ 
to the south pole of $S^2$. Since $\Delta^2$ deformation retracts onto the 
union of these two sides,  the mapping space can be retracted to the 
constant map $\Delta^2 \mapsto \ast$. 

We finish the argument by bootstrapping: the contractibility of
$|Y_\bullet|$ means that we may replace $K^2_\bullet$ by its fiber over
$Y_\bullet$.   The next step is to project this fiber to one of the other factors
of $S^2$, say with $v_1$. The corealization of this projection is also
contractible, since, unraveling definitions, it is the space of maps
from $\Delta^2$ to $S^2$ sending the $\{2\}$ edge (by definition) and the
$\{1\}$  edge (since we are in the fiber over $Y_\bullet$)
to the south pole, which is contractible.   
Repeating this with the $v_2$ projection onto $S^2$ we
see that $|K^2_\bullet|$ is indeed contractible.
\end{proof}

In our bootstrapping argument, it {\em does} matter onto 
which factor of $S^2$ we project $K^2_\bullet$ first.  
The image under the $v_1$ or $v_2$ projection does 
not have a corealization which is readily seen to be contractible 
without first looking in the fiber of the $\pi_{12}$ projection.  

Before proceeding to the next case, we must treat one last technical 
detail about \mbox{$\Delta^n$-spaces}, namely that of replacing a 
map by a fibration, as is standard for spaces.  

\begin{definition}\label{D:fibrep}
Let $f \colon X_\bullet \to Y_\bullet$ be a map of 
\mbox{$\Delta_n$-spaces} and define
the spaces and maps $X_\bullet \overset{i}{\to} \wt{X}_\bullet
\overset{\tilde{f}}{\to} Y_\bullet$ as follows.  Let $\wt{X}_S$ be 
the space of pairs $(x, \gamma)$ where $x \in X_S$ and 
$\gamma$ is a path in $Y_S$ with $\gamma(0) = f_S(x)$.  Inclusions 
between the $\wt{X}_S$ needed for the  $\Delta_n$-structure are 
defined by using those for $X_S$ and $Y_S$.  The map $i$ is 
defined by sending $x$ to $(x, c)$,  where $c$ is the constant path 
at $f_S(x)$.   The map $\wt{f}$ is defined by sending $(x, \gamma)$ 
to $\gamma(1)$.
\end{definition}

One may check that sending $(x, \gamma) \in \wt{X}_S$ to 
$x \in X_S$ defines a homotopy inverse to $i$ as 
\mbox{$\Delta^n$-spaces} so that $i$
induces an equivalence on corealizations.  
Moreover, it is easy to check that $\tilde{f}$ is a fibration.

\begin{theorem}\label{T:amthree}
The space $AM_3(\I^3)$ is homotopy equivalent to $\Omega^3 F$,
the space of based maps from $S^3$ to $F$, the homotopy fiber of the 
inclusion $S^2 \vee S^2 \hookrightarrow S^2 \times S^2$.  
\end{theorem}

Moreover, $F$ is known to be $\Sigma(\Omega S^2 \wedge \Omega S^2)$,
which splits as a countable wedge of spheres by the James splitting
(compare with Corollary~5.6 of \cite{GW99a}).

\begin{corollary}
The components of $AM_3(\I^3)$ are in bijective correspondence with the
integers.
\end{corollary}

\begin{proof}  
By the theorem, the components of $AM_3(\I^3)$ correspond to elements of
$\pi_3(F)$.  Since $\pi_3(S^2) = \Z$, we have 
$\pi_3(S^2 \times S^2) = \Z^2$ and 
by the Hilton-Milnor theorem \cite{Hi55}, 
$\pi_3(S^2 \vee S^2) = \Z^3$.  Because the
composites $S^2 \vee S^2 \hookrightarrow S^2 \times S^2 \to S^2$ are split, 
we see that $\pi_3(F) = \Z$.
\end{proof}

\begin{proof}[Proof of Theorem~\ref{T:amthree}]
We will be a bit more terse than in the previous proof since
the outline of the argument is precisely the same.

Once again we begin by defining a relatively simple $K^3_\bullet$ to which
$C^3_\bullet$ maps.
Set $K^3_{\emptyset} = (S^2)^6$ and define a map
from $C^3_\bullet$ over $K^3_\bullet$ by composing 
$\pi_{12} \times \pi_{23} \times \pi_{13} \times v_1 \times v_2 \times v_3$
with the collapsing map $\sigma$ (from the previous proof) on the first three
factors. We can give $K^3_\bullet$ the structure of a \mbox{$\Delta_3$-space} 
simply by defining the subspaces $K^3_S$ to be the images of the 
corresponding spaces in $C^3_\bullet$.   

The first step is to show that $|K^3_\bullet|$ is contractible.
Define a \mbox{$\Delta_3$-space} $Y_\bullet$ exactly as in the 
previous proof, by projecting $K^3_\bullet$ onto the first factor of 
$S^2$ (once again, the image of $\pi_{12}$).   The corealization
is manifestly contractible; since $Y_{\{2,3\}} = Y_{\{0\}} = \ast$,
we may define the contraction by deforming
$\Delta_3$ onto the union of this edge and face.
Thus we may consider the fiber of $K^3_\bullet$ over $Y_\bullet$.    Next
consider  the projections to $S^2$ defined by $v_1$ and $v_2$ in turn.
We claim the corealization of each is contractible; for instance, under
the projection $v_1$, the $\{0\}$ face (by definition) and 
$\{1\}$ face (by the bootstrap) map to $\ast$.   
The claim follows since $\Delta_3$ can be contracted 
to the union of these faces.
Proceeding in a symmetric way for the other factors of $S^2$
shows that $|K^3_\bullet|$ is contractible, as desired.

By construction, there is a projection map $\pi \colon C^3_\bullet
\to K^3_\bullet$ so in principal we should be able to identify $C^3_\bullet$ 
with the fiber of this projection.  But, this projection is not a fibration of
\mbox{$\Delta_3$-spaces}, since in particular it is not a fibration on 
$C^3_{\emptyset}$.  
Using Definition~\ref{D:fibrep} we may replace $C^3_\bullet$ by a
$\wt{C}^3_\bullet$ which has an equivalent corealization and which does 
fiber over $K^3_\bullet$ through a map $\tilde{\pi}$.  Investigating the fiber
$F_\bullet$ of $\tilde{\pi}$ we see that $F_S$ is contractible if 
$S$ is non-empty, since in these cases $\pi_S$ is a homotopy equivalence 
which implies  that $\tilde{\pi}_S$ is as well.  
A  \mbox{$\Delta_n$-space} $X_\bullet$  with $X_S$ contractible for all 
non-empty $S$ has a corealization which is homotopy equivalent to 
$\Omega^n X_{\emptyset}$.  
Thus, the corealization of $F_\bullet$ is $\Omega^3 F_{\emptyset}$.  
But by definition $F_{\emptyset}$ is the homotopy fiber of the projection 
$\pi_{12}  \times \pi_{23} \times \pi_{13}$ from $C_3(\R^3)$ to $(S^2)^3$.  
By looking at the fiber of $\pi_{12}$ (which is by itself a fibration) this is 
the same as the fiber  of $\pi_{13} \times \pi_{23}$ from  
$\R^3 \setminus \{a_1,a_2\}$ to $S^2 \times S^2$ where in this case
$\pi_{13}$ sends $z \in \R^3 \setminus \{a_1,a_2\}$ to the unit vector 
from $a_1$ to $z$ (and similarly for $\pi_{23}$).  Retracting 
$\R^3 \setminus \{a_1,a_2\}$ onto $S^2 \vee S^2$ 
we see that $\pi_{13} \times \pi_{23}$ coincides with the inclusion  
$S^2 \vee S^2 \hookrightarrow S^2 \times S^2$. 
\end{proof}

\section{Linking of \colinearity\ manifolds}\label{definition}

\subsection{Definition of $\nu_2$}
In this section we formalize the definition of our self-linking invariant 
$\nu_2$.  
In order to proceed we need to 
address technicalities about collinearity submanifolds of compactified 
configuration spaces.

\begin{definition}
Define $\coll_i (\I^k)$ to be the subspace of $C_3(\I^k)$ consisting 
of configurations of points  $x=(x_1,x_2,x_3)$ such that  
$\{x_1,x_2,x_3\} \subseteq L$ where $L$ is a  straight line in  
$\R^k$, and in that line $x_i$ is between the other two points.  
Alternatively,   $\coll_i (\I^k)$ is the preimage of the submanifold of 
$S^{k-1} \times S^{k-1}$ defined   by pairs $u \times -u$ under the 
map $\pi_{ij} \times \pi_{i\ell}$.
\end{definition}

As mentioned in Theorem~\ref{T:compact}, the maps $\pi_{ij}$ extend to 
smooth maps from $C_n[\I^k, \partial]$ to $S^{k-1}$ simply by restricting 
the projection of $B_{n,k}$ onto the $ij$ factor of $S^{k-1}$ .  
Let $\wt{\pi}_{ij} : \wt{C}_n(\R^k) \to S^{k-1}$ be the unique map
which factors $\pi_{ij}$ through $\wt{C}_n(\R^k)$.  
Using the decomposition of Theorem~\ref{T:strata} and the coordinates 
around these strata given in \cite{Si02a},
we deduce that  $\pi_{ij} \times \pi_{ik}$ restricted to a stratum 
$S_\mathcal{P}$  on the boundary of $C_3[\I^k]$ is either projection 
onto a factor of  $\wt{C}_2(\R^k) \times \wt{C}_2(\R^k)$ or factors as 
projection onto  $\wt{C}_3(\R^k)$ followed by 
$\wt{\pi}_{12} \times \wt{\pi}_{13}$.  
We deduce that the maps $\pi_{ij}$ are open maps.  
We obtain the following proposition as a consequence of transversality.

\begin{proposition}
The submanifolds $\coll_i (\I^k)$ extend to be submanifolds with
corners of $C_n[\I^k, \partial]$, which we denote $\coll_i[\I^k]$.
\end{proposition}

It will be helpful at times to be more explicit about some of the boundary of
$\coll_i[\I^k]$.   

\begin{proposition}\label{collinearboundary}
The submanifold $\coll_i[\I^k]$ has the following intersections with boundary strata.
\begin{itemize}
\item In $S_{(i=j)}$, consisting of pairs 
$(x_1, x_2)  \times ((y_1, y_2)/\sim) \in C_2(\I^k) \times \wt{C}_2(\R^k)$ 
such that  $x_2 - x_1 = -k (y_2 - y_1)$ where $k>0$.
 \item In $S_{(1=2=3)}$, consisting of equivalence classes in 
$\wt{C}_3(\R^k)$  which are represented by collinear triples.
\item In $S_{((i=j)=k)}$, consisting of points represented by 
$x, v, -v$ in  $\I^3 \times \wt{C}_2(\R^k) \times \wt{C}_2(\R^k)$.
\end{itemize}
\end{proposition}

For the sake of transversality arguments we choose the boundary points 
on $\I^3$, as needed for the definitions for
${\rm Emb}(\I, \I^3)$ and $AM_n(\I^3)$, to be $(\frac{1}{2}, \frac{1}{2},0)$ and 
$(\frac{1}{2} + \epsilon, \frac{1}{2}, 1)$ for some $\epsilon > 0$.
We choose the tangent vectors at those points to be perpendicular 
to the boundary.

\begin{theorem}\label{trans1}
Any $g \in AM_3(\I^3)$ is arbitrarily close to some $g' \in AM_3(\I^3)$ 
which is  transverse to $\coll_i[\I^3]$ for $i = 1, 3$.
\end{theorem}

\begin{proof}
We may establish this inductively over the strata of $A_3$ using
the Extension Theorem from \cite{GP74}.  The images under $g$ of 
the vertices of $A_3$,  which are fixed by choosing boundary points 
and vectors in $\I^3$, are disjoint from $\coll_i[\I^3]$ by our choice above.  
Inductively, we may deform $g$ slightly to be  transverse to
$\coll_i[\I^3]$ on the dimension $i$ strata, having fixed it on the 
dimension $i-1$ strata, by using Proposition~\ref{collinearboundary}
to check that the intersections of 
$\coll_i[\I^3]$ with the aligned strata of $C_3[\I^3, \partial]$ have 
codimension two, which is true in all cases except when $i = 2$ 
and the stratum in question is labeled by $(1=2=3)$ or any 
parenthesization which includes $(1=2=3)$.
\end{proof}

\begin{definition}
Given $g \in AM_3(\I^3)$, deform $g$ to be transverse to  
$\coll_i[\I^3]$ for $i = 1,3$  and define $\coll_i[g]$, the subspaces
of collinear points, to be $g^{-1}(\coll_i[\I^3]) \subset A_3$.
\end{definition}

In the interest of defining knot invariants directly, we also show 
that the evaluation  map of a knot can be made transverse to 
collinearity submanifolds by a small isotopy of  the knot.

\begin{theorem}
For a generically positioned knot $f$,  $C_3[f]$ is transverse to $\coll_i [\I^3]$ 
for $i\in \{1,3\}$.
\end{theorem}

\begin{proof}
Let $sk_i(A_n)$ be the $i$-skeleton of $A_n$.  For any knot $f : \I \to \I^3$,
 $C_3[f]|_{sk_0(A_3)}$ is disjoint  from $\coll_i[\I^3]$ by our choice of 
endpoint data for $f$.   Transversality of $C_3[f]$ over 
$sk_1(A_3)$ can be guaranteed by   demanding that the tangent lines to 
$f$ do not intersect $f(0)$ or $f(1)$,  a generic condition.

For the two-skeleton investigate the transversality face-by-face.
Collinear points in  $A^{(0=1)}_3$ correspond to collinearities on the knot 
which include $f(0)$, so one is guaranteed transversality if the
projection of the knot to the sphere based at $f(0)$ is a regular knot diagram,
 which is well-known to be a generic condition.  Transversality in 
$A^{(3=4)}_3$ works similarly.

Collinearities in the $A^{(1=2)}_3$ and $A^{(2=3)}_3$ correspond, 
combining Proposition~\ref{collinearboundary} with the definition 
of the evaluation  map, to pairs of times $s, t$ such that the tangent 
line of $f$ at $s$ intersects the knot at $f(t)$.  
Informally, transversality requires that small variations of $s$ and $t$ 
should give  rise to a two-dimensional family of non-collinear triples.  
Formally, this  means that the three vectors $f(s)-f(t)$, $f'(t)$ and $f''(s)$ 
are  linearly independent.
This condition has a simple geometric interpretation. 
Provided the knot $f$  has everywhere non-zero curvature, the map 
$\exp : \I \times (\R-\{0\}) \to \R^3$ given by
$\exp(t,h) = f(t)+hf'(t)$ is an immersed submanifold of $\R^3$. The above
transversality condition is equivalent to the condition that $f$ intersects
the image of $\exp$ transversally, which is a generic condition. The fact
that non-zero curvature is a generic condition for knots follows quickly
from the Frenet-Serret Theorem \cite{MP77}.

Finally, for an interior point $(t_1,t_2,t_3) \in A^\phi_3$ the transversality 
condition is that for a \colinear\ triple, the four vectors $f(t_1)-f(t_2)$,
$f'(t_1)$,$f'(t_2)$,$f'(t_3)$ must
span $\I^3$.  While difficult to see intuitively, this is 
straightforward to prove formally.
Consider the map $\Phi : SO_3 \times \I^3 \to (\I^2)^3$  which sends
a rotation matrix $A \in SO_3$ and three times $t_1,t_2,t_3$ to
the orthogonal projections of 
$A(f(t_1)),A(f(t_2)),A(f(t_3))$ to $(\I^2\times\{0\})^3 \subseteq (\I^3)^3$. 
Our transversality condition follows if $\Phi$ transversally intersects the 
diagonal $x_1=x_2=x_3$ in $(\I^2\times\{0\})^3$. Apply a standard
transversality argument, extending $\Phi$ to a function 
$\overline{\Phi} : SO_3 \times \I^3 \times \I^k \to (\I^2\times \{0\})^3$, 
such that $\overline{\Phi}_{|SO_3 \times \I^3 \times \{0\}} = \Phi$,
which is transverse to the diagonal.  The Transversality Theorem 
\cite{GP74} tells us that an arbitrarily small isotopy of $f$ satisfies 
our transversality condition.
\end{proof}

Finally, we analyze where the boundaries of these collinearity submanifolds may lie.

\begin{proposition}
Any $g \in AM_3(\I^3)$ may be deformed slightly so that the $1$-manifold
$\coll_1[g]$ has boundary only on $A^{(1=2)}_3$ and $A^{(3=4)}_3$.  
Similarly, $\coll_3[g]$ may be assumed to have boundary only on 
$A^{(0=1)}_3$ and $A^{(2=3)}_3$.
\end{proposition}

\begin{proof}
We focus on $\coll_1[g]$ since the arguments for $\coll_3[g]$ are 
identical up to re-indexing and ``turning $\I^3$ upside-down.''  We
apply Proposition~\ref{collinearboundary} and analyze only the codimension
 one faces of $A_3$ since by Theorem~\ref{trans1} there are no boundary 
 points on strata with higher codimension.

To set notation, let $a \in A_3$ and  
let $p(g(a)) = (x_1, x_2, x_3) \in (\I^3)^3$.
By definition of an aligned map, the images of the $A^{(0=\cdots=3)}_3$,
$A^{(1=\cdots=4)}_3$ and $A^{(1=2=3)}_3$ under $g$ 
are all in $\coll_2[\I^3]$, so $\coll_1[g]$ has no boundary
on these strata.  Next, $g(a)$ for $a \in A^{(0=1=2)}_3$ has
$x_1 = x_2 = \left(\frac{1}{2}, \frac{1}{2}, 0\right)$ and 
$\pi_{12}$ pointing downwards.
Thus, a $\coll_1[g]$-collinearity would force $g(t_3)$ to be above $\I^3$,
which of course is not allowed.  
On both $A^{(2=3=4)}_3$ and $A^{(2=3)}_3$, 
we have $\pi_{12}(g(a)) = \pi_{13}(g(a))$, which means that 
these unit vectors  cannot be opposites as required for $\coll_1[g]$.  
Finally, on $A^{(0=1)}_3$ we have  $x_1(g(a)) = 
\left( \frac{1}{2}, \frac{1}{2}, 0 \right)$ 
so that a $\coll_1[g]$-collinearity
would require that  both $x_2$ and $x_3$ have last coordinate
zero and be on  a line, which is a codimension three condition.
\end{proof}

Because  $\coll_1[f] $ and $\coll_3[f]$ have no boundary on the
hidden faces of $C_3[\I, \partial]$, they project to
manifolds with boundary in $\Delta^3$. By abuse, we give this projection the same 
name, and for the rest of the paper we will work in the simplex $\Delta^3$. 
Moreover $\coll_i (\R^k)$, and thus $\coll_i (\I^k)$, is orientable 
since it is diffeomorphic to the product 
$(\R^k) \times S^{k-1} \times (0,\infty)^2$.
Thus $\coll_i[\I^k]$ is orientable,
and we fix such an orientation for the remainder of this discussion.  
By pulling back the orientation of the normal bundle of $\coll_i [\I^3]$ in 
$C_3[\I^3, \partial]$, we get an orientation of $\coll_i[g]$ for any 
$g \in AM_3(\I^3)$.

Instead of labeling the faces of $\Delta^3$ by subsets of $\{0, \ldots, 3\}$,
as in the definition of $\Delta^3_\bullet$, we now use parenthesizations, 
as inherited from  the projection from $A_3$.

\begin{definition}
Define a closure of the manifold with boundary $\coll_1[g]$ 
to be any closed piecewise smooth $1$-manifold $\overline{\coll_1[g]}$ 
whose intersection with ${\rm Int}(\Delta^3)$ is $\coll_1(g)$, and whose 
intersection  with $\partial(\Delta^3)$ lies only in $\Delta^3_{(1=2)}$ 
and $\Delta^3_{(3=4)}$. Define $\overline{\coll_3[g]}$ similarly. 
\end{definition}

A Seifert surface for
$\overline{\coll_1[g]}$ can only intersect $\overline{\coll_3[g]}$ inside the
simplex, and these completions have an orientation determined by the
orientations on $\coll_i[g]$.  Thus,  the linking number of these  two manifolds
is defined and will be  independent of choices of these completions. 
We may now define our invariant.

\begin{definition}
For $g \in AM_3(\I^3)$, let $\mu_2(g) =
lk(\overline{\coll_1[g]},\overline{\coll_3[g]})$.  For a knot, define
$\nu_2(\knot)$ to be $\mu_2(C_3[\knot])$.
\end{definition}

A homotopy between $f$ and $g$ in $AM_3$ may again be deformed 
to be transverse to the \colinearity\ conditions, and thus give rise to 
an oriented cobordism  between the manifolds
$\coll_i[f]$ and $\coll_i[g]$.  Thus, $\mu_2$ passes to a map 
$\pi_0(AM_3(\I^3)) \to \Z$.  We next establish that 
this invariant encodes all of the  information we may get about
knotting in $\I^3$ from the evaluation map $C_3[f]$.

\subsection{$\mu_2$ is an isomorphism}

By \refT{amthree}, $\pi_0(AM_3(\I^3)) \cong \pi_3(F)$, 
where $F$ is the homotopy fiber of the inclusion from
$S^2 \vee S^2$ to $S^2 \times S^2$.  The composites
$S^2 \vee S^2 \to S^2 \times S^2 \to S^2$, where the second
map is projection onto a factor, are split so that $\pi_3(F)$
is a subgroup of $\pi_3(S^2 \vee S^2)$.  That subgroup is
isomorphic to the integers, with the isomorphism realized by
sending an $f \colon S^3 \to S^2 \vee S^2$  to
$lk_{S^3}\left(f^{-1}(a),f^{-1}(b)\right)$, where
 $a\in (S^2\vee \ast)\backslash \{\ast\}$ and $b\in (\ast\vee S^2)\backslash
\{\ast\}$ are regular values of $f$.
This observation led us to our definition of $\mu_2$, and 
one could trace through the equivalence given by 
\refT{amthree} in order to equate $\mu_2$ with this linking number
and thus show that $\mu_2$ is an isomorphism.  
To do this a good choice of embedding 
$S^2\vee S^2\subset \R^3\setminus\{a,b\}$ is with one sphere is centered at $b$ of radius 
$|a-b|/2$, and the second  sphere is centered at $a$ of radius $3|a-b|/2$.
We prefer an indirect
but shorter approach, which takes advantage of our computations
from Section~\ref{examples}.

\begin{definition}
Fix a generic embedding of the unknot $U$ in $\I^3$, whose
collinearity submanifolds are empty.
Let $\phi|C^3_\bullet|$ and $\phi|F^3_\bullet|$ denote the
subspaces of $|C^3_\bullet|$ and $|F^3_\bullet|$ respectively
of maps whose restriction to $\partial \Delta^3$  coincides with 
the restriction to $\partial \Delta^3$ of $C_3[U]$, 
or respectively coincides on $\partial
\Delta^3$ with a chosen lift of $C_3[U]$ to $|F^3_\bullet|$.
\end{definition}

\begin{theorem}\label{T:comd}
There is a commutative diagram 
\begin{equation*}
\begin{CD} 
\pi_0(|F^3_\bullet|) = \Z  @>{(1)}>> \pi_0(|C^3_\bullet|) @>{\mu_2}>> \Z \\
@A{(2)}AA @AAA @A{(5)}AA\\
\pi_0(\phi|F^3_\bullet|) @>{(4)}>> \pi_0(\phi |C^3_\bullet|) @>{(3)}>> \pi_3(C_3\la \I^3 \ra)
\end{CD}
\end{equation*}
where the maps (1) and (2) are isomorphisms of sets and the maps 
(3), (4), and (5) are homomorphisms.
\end{theorem}

We first indicate how this theorem leads to the main result
of this subsection.

\begin{corollary}
The map $\mu_2$ is an isomorphism.
\end{corollary}

\begin{proof}
Chasing through the diagram in \refT{comd} we see that 
$\mu_2$ is a homomorphism since it is the composite
of two bijections, the inverses to (1) and (2), with the
homomorphisms (3), (4), and (5).

In \refT{comp} below, we show that $\mu_2$ takes on the
value one (for the evaluation map of the trefoil knot)
and thus must be an isomorphism.
\end{proof}

\begin{proof}[Proof of \refT{comd}] 

We define the maps in the diagram.

\begin{enumerate}

\item The map (1) is the map on $\pi_0$ of the induced map on 
corealizations of the composite map of $\Delta_n$-spaces 
$F_\bullet \to \wt{C}^3_\bullet \to C^3_\bullet$, which is
shown to be an equivalence in \refT{amthree}.

\item Because all of the boundary terms of $F^3_\bullet$
are contractible, we may induct over the skeleta of
$\partial \Delta^3$ to show that the inclusion of
$\phi|F^3_\bullet|$ in $|F^3_\bullet|$ is an equivalence,
and thus define (2) as the induced isomorphism.

\item The map (3) is a case of a general construction which we now
recall.  

If $\mathcal{M}$ is a space of maps from $\Delta^n$
to $X$ whose restriction to $\partial \Delta^n$ is prescribed,
$\mathcal{M}$ is homotopy equivalent to $\Omega^n X$ as follows.
Fix one particular $f \in \mathcal{M}$.  Define the ``gluing
with $f$'' map $\Gamma_f  : \mathcal{M} \to \Omega^n X$ 
by having $\Gamma_f(g) : S^n = \Delta^n \cup_\partial \Delta^n \to X$
restrict to $f$ on the first $\Delta^n$ and $g$ on the second.

The map (3) is $\pi_0 (\Gamma_{C_3[U]})$.  It is an isomorphism of
sets since $\Gamma_{C_3[U]}$ is a homotopy equivalence.  In fact
we use this isomorphism to  define the group structure on 
$\pi_0(\phi|C^3_\bullet|)$.

\item The map underlying (4) is simply the restriction of the map 
underlying (1) to $\phi|F^3_\bullet|$, which maps to 
$\phi|C^3_\bullet|$ by construction.  It induces a homomorphism
on $\pi_0$ since it sits in a commutative square
\begin{equation*}
\begin{CD} 
\phi|F^3_\bullet| @>{(4)}>> \phi |C^3_\bullet| \\
@V{\Gamma_{\wt{C_3[U]}}}VV @V{\Gamma_{C_3[U]}}VV \\
\Omega^3 F_\phi @>>> \Omega_3(C_3\la \I^3 \ra),
\end{CD}
\end{equation*}
where the vertical arrows are used to define the group structures on
$\pi_0$ of $\phi|F^3_\bullet|$ and $\phi|C^3_\bullet|$ and the bottom
arrow is a homomorphism on $\pi_0$ since it is a map of three-fold
loop spaces.

\item Finally, the map (5) is defined by taking an 
$f : S^3 \to C_3 \la \I^3, \partial \ra$, homotoping it so that
its image is in $C_3 (\I^3)$ (whose inclusion in $C_3 \la \I^3, \partial \ra$
is a homotopy equivalence) and transverse to $Co_1$ and $Co_3$, and then taking
$lk_{S^3}(Co_1(f), Co_3(f))$. 
In general, sending $\pi_n(M)$ to $\Z$ by taking the linking number of the
preimages of two disjoint closed submanifolds of $M$ (which itself need not
be closed) whose codimensions add to $n+1$ is a homomorphism, so in
particular (5) is a homomorphism.
\end{enumerate}

It is straightforward to check that the maps as constructed commute
since $C_3[U]$ has empty collinearity submanifolds,
and this completes the proof.

\end{proof}

\section{Examples}\label{examples}

Clearly if $\knot$ is the unknot, $\nu_2(\knot)=0$,
since any reasonably simple generic embedding of the unknot has no 
\colinear\ triples.
In this section we will compute $\nu_2(\knot)$ directly for the trefoil 
and figure eight knots.   To do so, one could parametrize these knots and solve 
the systems of equations which arise from collinearity conditions.  Indeed 
as explained further at the end of Section~\ref{coquad}, one 
application of our work is to, for example, bound the value of 
the degree two Vassiliev invariant of a knot which is parametrized by polynomials 
of a given degree.  We prefer to take a 
slightly less explicit but more geometric approach by choosing embeddings with 
certain monotonicity properties, which are depicted in figures \ref{trefoil}  
and \ref{fig8}.  Throughout this section we give full arguments for the
trefoil but leave the entirely analogous arguments for the figure eight
to the reader.  This analysis will proceed first by finding boundary points of 
the collinearity submanifolds, and will ultimately highlight the fact that 
quadrisecants of the knot may play an essential role in computation
 of the invariant, a theme which is fully developed in the next section.

\begin{center}
\begin{minipage}{7cm}\begin{mydiagram}\label{trefoil}\begin{center}
{
\psfrag{f1/6}{$f(\frac{1}{6})$}
\psfrag{f1/3}{$f(\frac{1}{3})$}
\psfrag{f1/2}{$f(\frac{1}{2})$}
\psfrag{f2/3}{$f(\frac{2}{3})$}
\psfrag{f5/6}{$f(\frac{5}{6})$}
$\includegraphics[width=7cm]{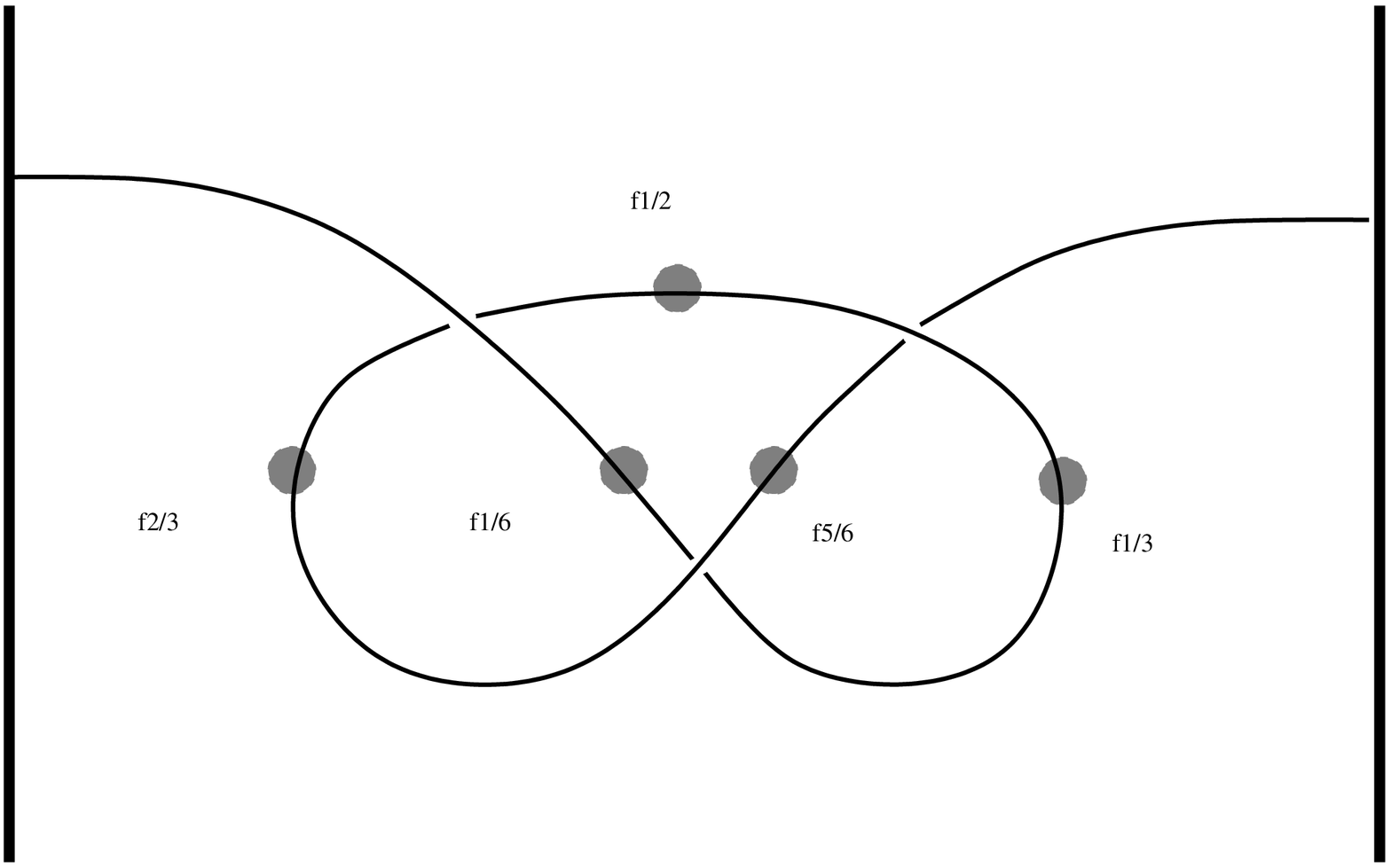}$
}
$x_3=0$ points on the trefoil knot
\end{center}\end{mydiagram}\end{minipage}
\begin{minipage}{7cm}\begin{mydiagram}\label{fig8}\begin{center}
{
\psfrag{f1/8}{$f(\frac{1}{8})$}
\psfrag{f1/4}{$f(\frac{1}{4})$}
\psfrag{f3/8}{$f(\frac{3}{8})$}
\psfrag{f1/2}{$f(\frac{1}{2})$}
\psfrag{f5/8}{$f(\frac{5}{8})$}
\psfrag{f3/4}{$f(\frac{3}{4})$}
\psfrag{f7/8}{$f(\frac{7}{8})$}
$\includegraphics[width=7cm]{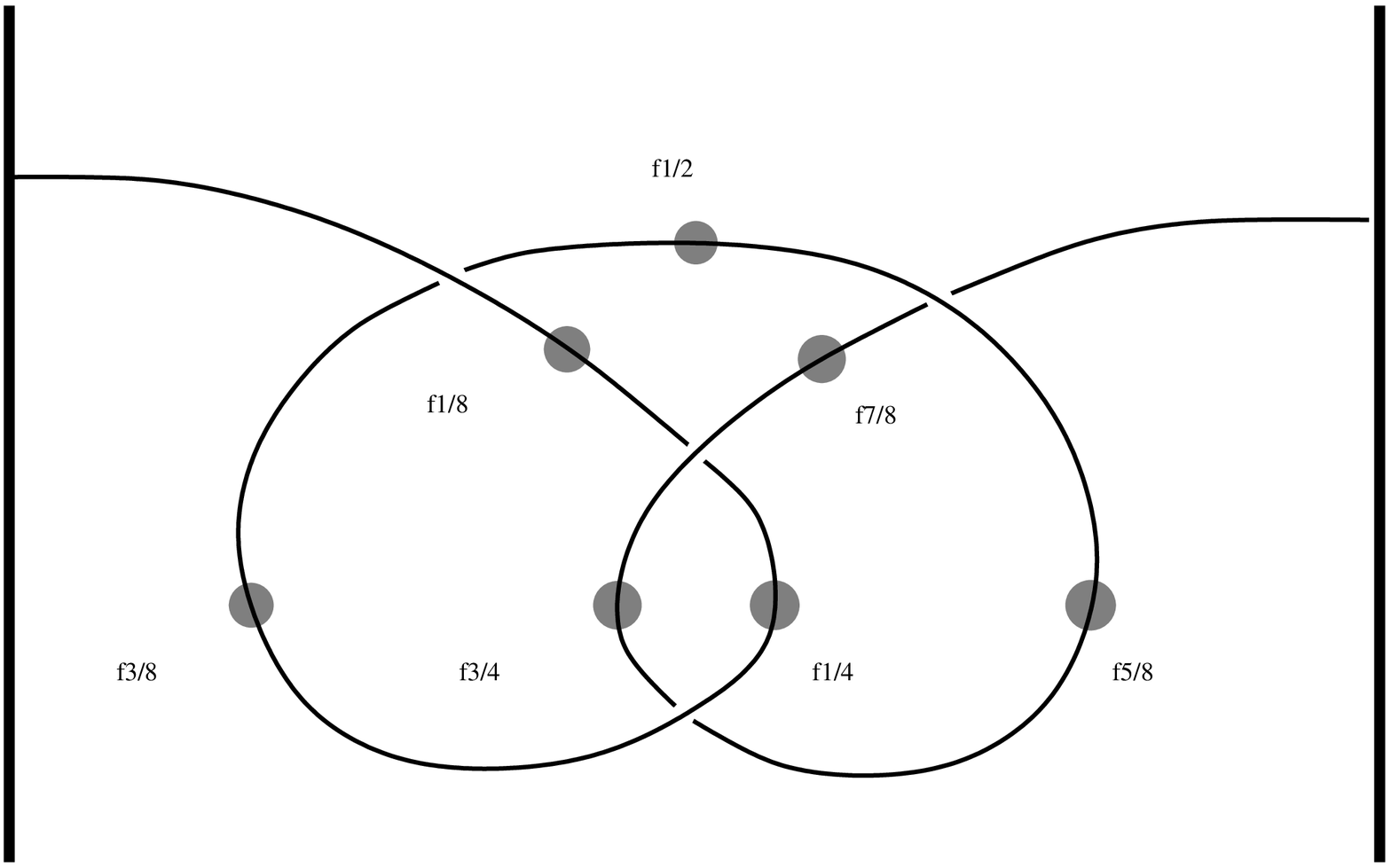}$
}
$x_3=0$ points on the figure-$8$ knot.
\end{center}\end{mydiagram}\end{minipage}\end{center}

The plane into which the knot is projected will be given the coordinates 
$x_1$ and $x_2$, and the coordinate that points out of the $(x_1,x_2)$-plane 
will be the $x_3$ coordinate.  In Figures \ref{trefoil} and \ref{fig8} 
the $x_3=0$ points of the embedding are specified by dots, five for the 
trefoil and  seven for the figure eight.  It will simplify the complexity of the 
collinearity submanifolds to assume that between each of the dots the
coordinate function  $x_3(t)$ strictly increases and then strictly 
decreases or vice-versa.

\begin{center}
\begin{minipage}{6cm}\begin{mydiagram}\label{trefoil_boundary_triples}
\begin{center}
$$\includegraphics[width=6cm]{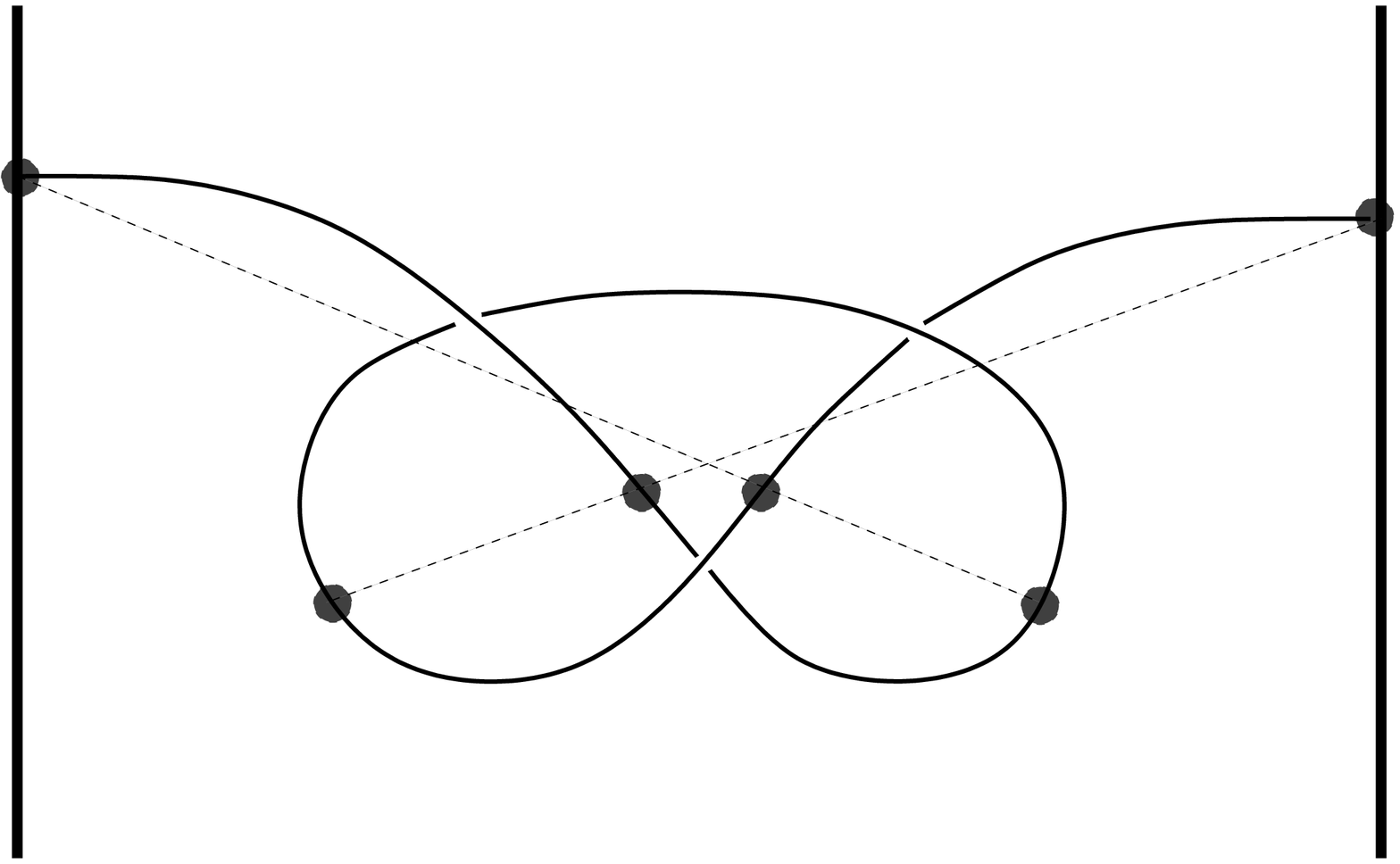}$$
\colinear\ triples on $(0=1)$ and $(3=4)$ strata.
\end{center}\end{mydiagram}\end{minipage}
\hskip 1.5cm
\begin{minipage}{6cm}\begin{mydiagram}\label{fig8_boundary_triples}
\begin{center}
$$\includegraphics[width=6cm]{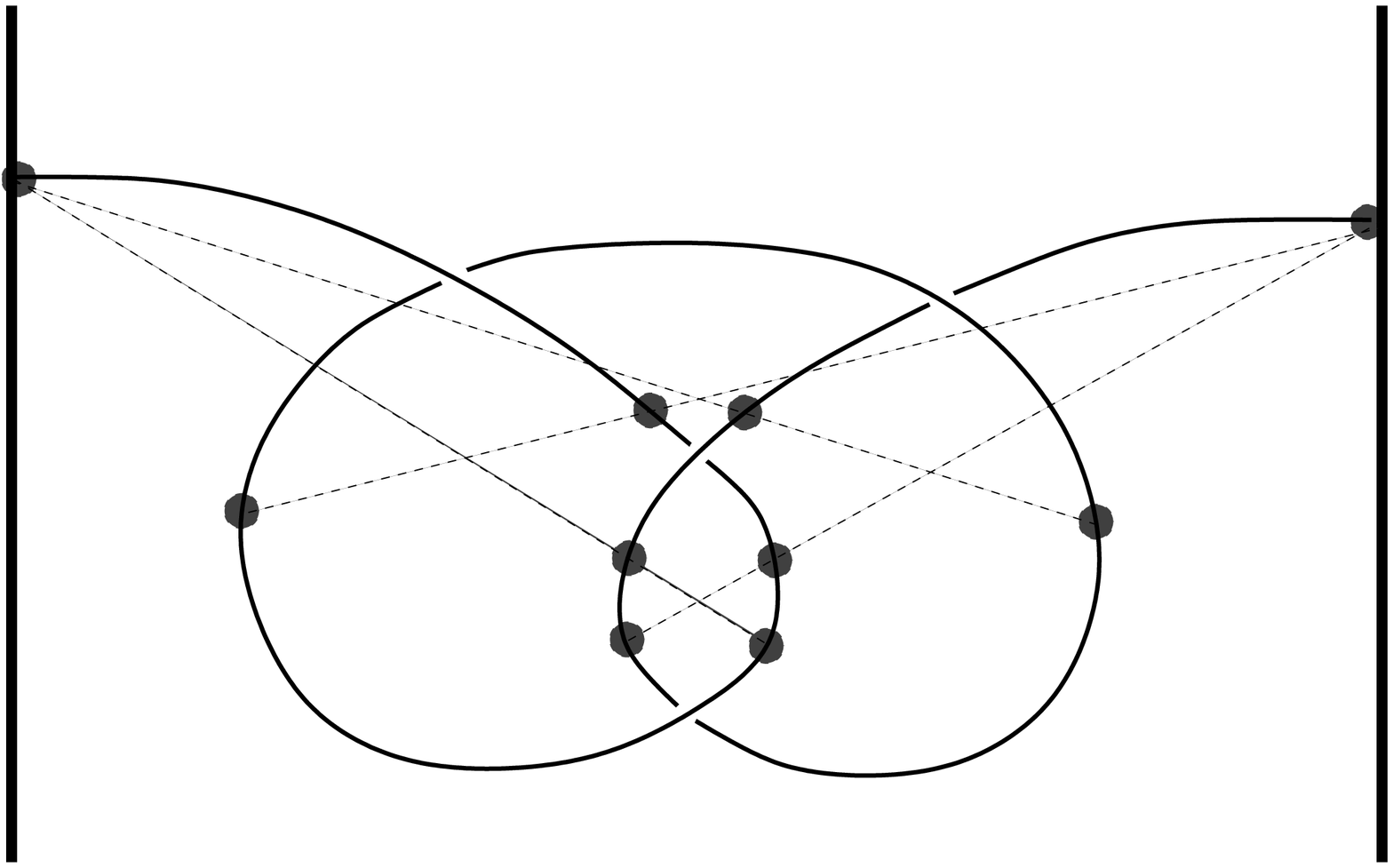}$$
\colinear\ triples on $(0=1)$ and $(3=4)$ strata.
\end{center}\end{mydiagram}\end{minipage}
\end{center}

In Figures \ref{trefoil_boundary_triples} and \ref{fig8_boundary_triples} 
all the \colinear\ triples $(f(t_1), f(t_2), f(t_3))$ where either $t_1=0$ or 
$t_3=1$ are indicated.  On $A^{(0=1)}_3$, one can label any \colinear\ 
triple by the number  $t \in \I$, such that the line segment
$[f(0),f(t)]$ contains one point of the knot $\knot$ in its interior.
In the case of the trefoil there are no \colinear\ triples
for $t \in [0,\frac{1}{6}]$ since $x_3(t)$ is decreasing on that interval. 
There is the one solution in $[\frac{1}{6},\frac{1}{3}]$ as sketched
above, since the knot near $t = \frac{5}{6}$ sits over the line from $f(0)$ to 
$f(\frac{1}{6})$ but under the line from $f(0)$ to $f(\frac{1}{3})$.  There 
are no solutions for $t \geq \frac{1}{3}$ except possibly one in $[\frac{5}{6},1]$, 
which would give rise to a point of $\coll_2 [f]$, and so is ignored.  

\begin{center}
\begin{minipage}{6cm}\begin{mydiagram}\label{trefoil_tangential_triples} 
\begin{center}
$$\includegraphics[width=6cm]{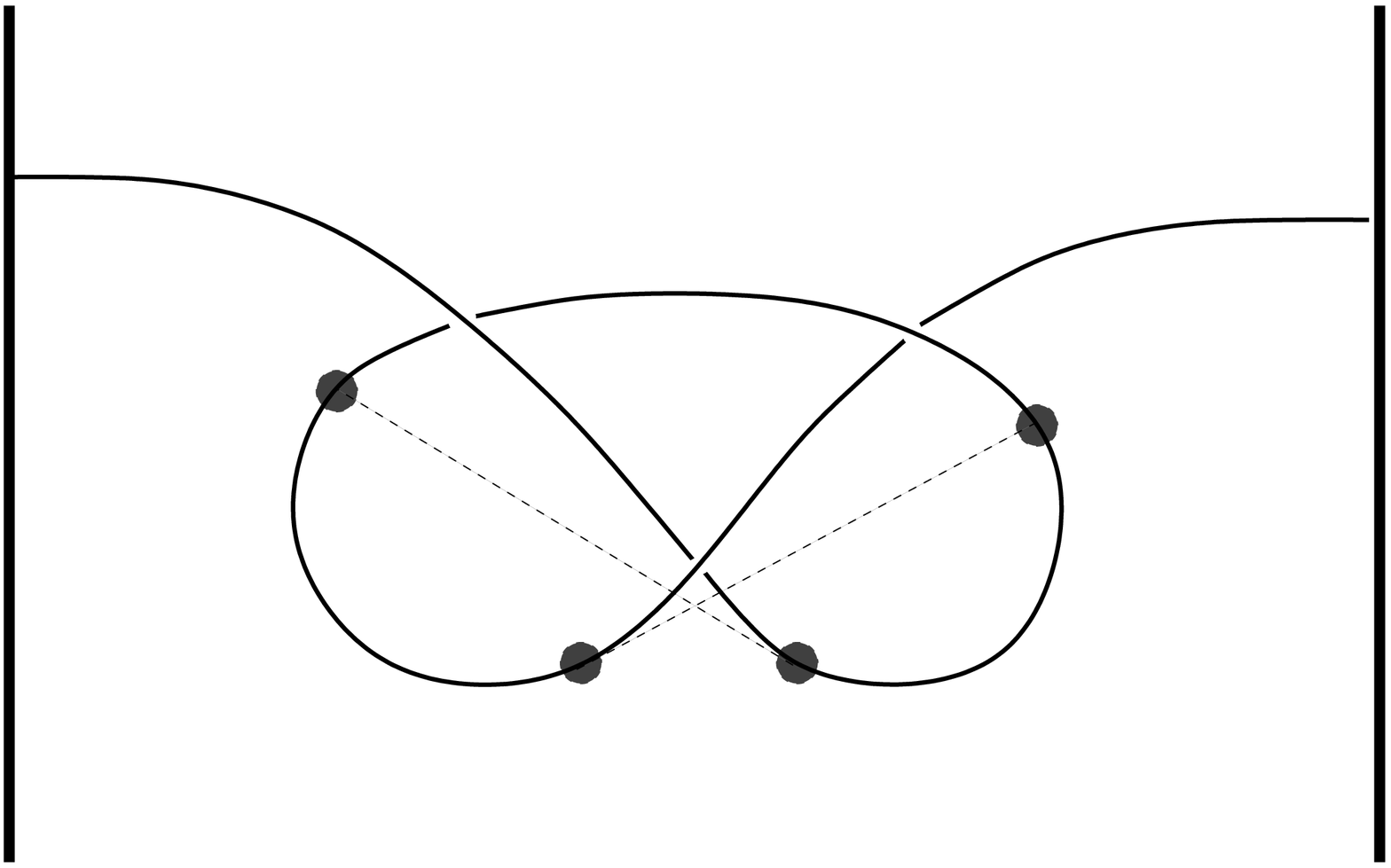}$$
\colinear\ triples on the $(1=2)$ and $(2=3)$ strata.
\end{center}\end{mydiagram}\end{minipage}
\hskip 1.5cm
\begin{minipage}{6cm}\begin{mydiagram}\label{fig8_tangential_triples}
\begin{center}
$$\includegraphics[width=6cm]{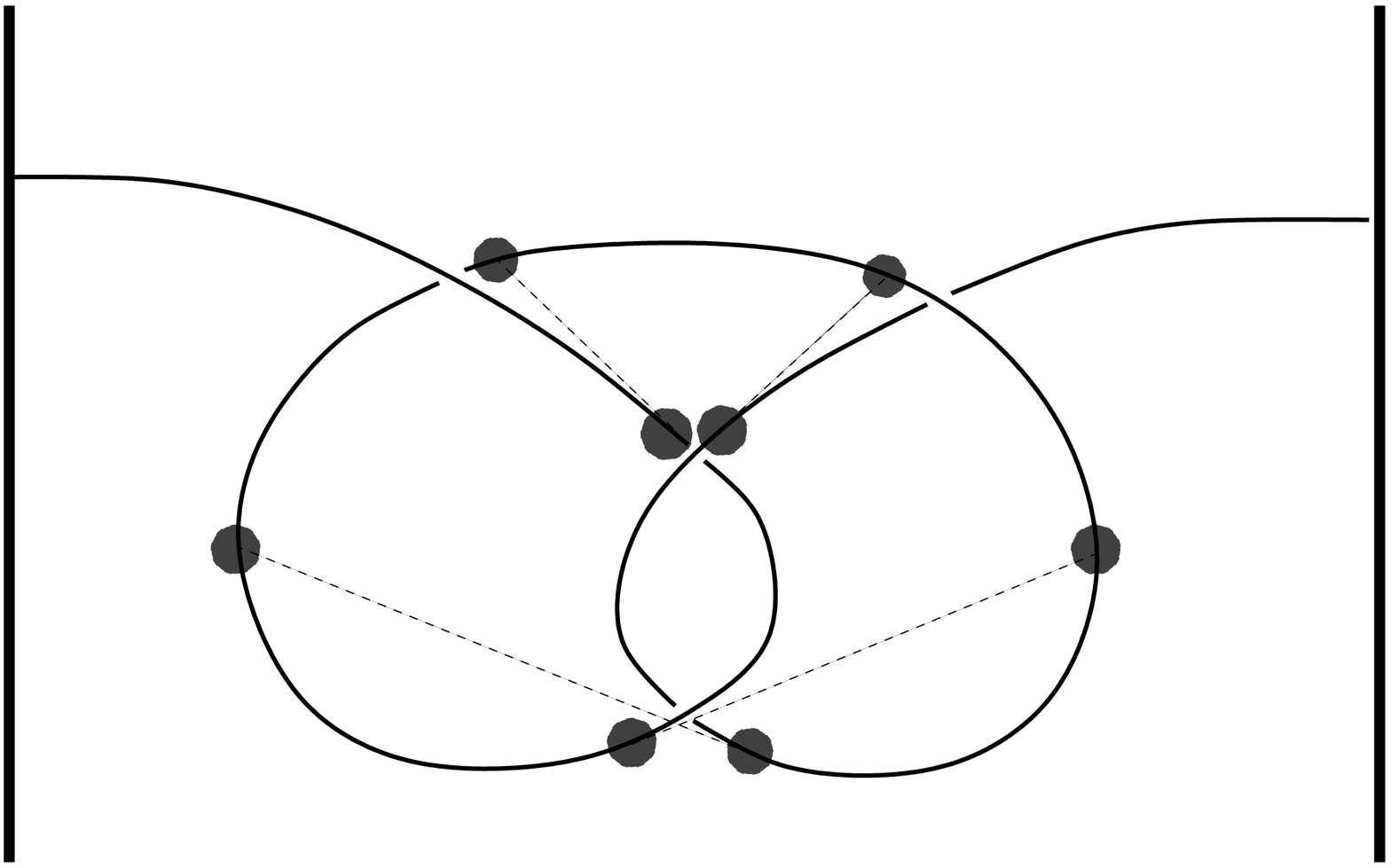}$$
\colinear\ triples on the $(1=2)$ and $(2=3)$ strata.
\end{center}\end{mydiagram}\end{minipage}
\end{center}

In Figures \ref{trefoil_tangential_triples} and  \ref{fig8_tangential_triples} all 
the tangential \colinear\ triples,  elements of $\coll_i[f]$ for $i\in \{1,3\}$ 
which lie in the $A_3^{(1=2)}$ and $A_3^{(2=3)}$ are indicated.   In all
these figures, points corresponding to  $\coll_2 [f]$ are systematically ignored.
On these strata collinear triples coincide with 
$t$ such that the tangent line to $\knot$ at $f(t)$ intersects the knot at a point 
other than $f(t)$.  In the 
case of the trefoil, we can deduce from boundary value and monotonicity 
arguments that there are tangential \colinear\ triples corresponding to
$t \in [\frac{1}{6},\frac{1}{3}]$ and $t \in [\frac{2}{3},\frac{5}{6}]$
as sketched above.  The tangential \colinear\ triples for 
$t \in [0,\frac{1}{6}]$ or $[\frac{5}{6}, 1]$,  would
belong to $\coll_2 [f]$ and can thus be disregarded. 
 
Now that we have the boundary structure of the $\coll_i[f]$ for the trefoil and
figure-eight knots, we need to understand the interior structure.  As usual, 
this can be understood from the crossing structure of a projection of these links.
The following lemma is simple and extremely useful.

\begin{lemma}\label{crossing}
Let $\rho \colon \Delta^3 \to \Delta^2$ be defined by forgetting the $t_1$ 
coordinate,  and let $f$ parametrize a knot $\knot$.
A crossing of the projection of $\coll_1[f]$ and $\coll_3[f]$ under $\rho$ 
corresponds  to a collinearity of four points on $\knot$.
\end{lemma}

\begin{proof}
A crossing of $\coll_3[f]$ and $\coll_1[f]$, corresponds to having
$f(t^*_1),f(t^*_2),f(t^*_3)$ on $L^*$ and  $f(t_1'),f(t_2'),f(t_3')$ on $L'$ with 
$t^*_2=t'_2$ and $t^*_3=t'_3$, since their projections under $\rho$ agree.
Because  $t^*_2=t'_2$ and $t^*_3=t'_3$, we have $L^* = L'$ so in fact the four points 
$f(t_1'), f(t_1^*), f(t_2')$ and $f(t_3')$ are all collinear.
\end{proof}

\begin{center}
\begin{minipage}{6cm}\begin{mydiagram}\label{trefoil_quadruples}
\begin{center}
$$\includegraphics[width=6cm]{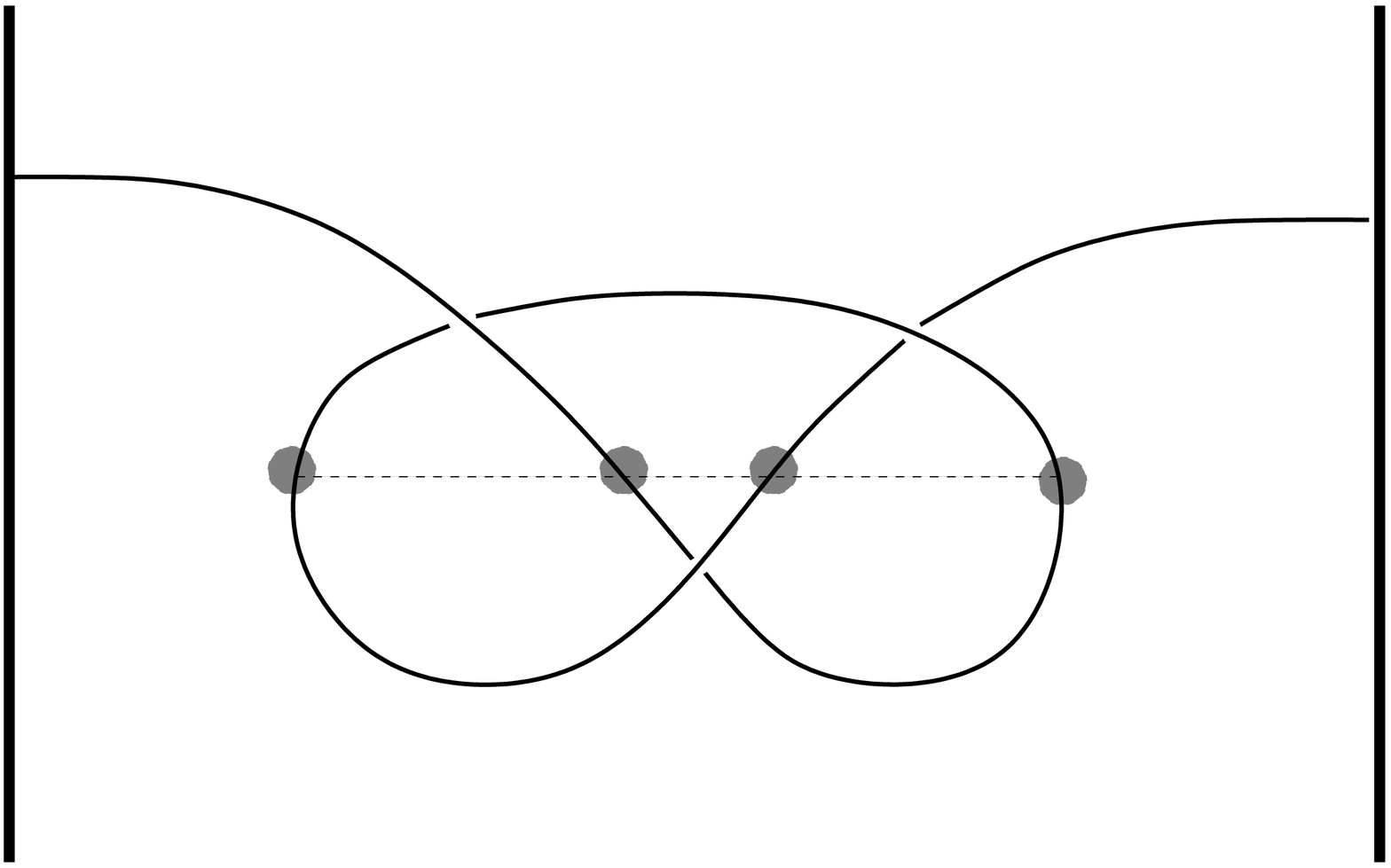}$$
quadrisecants.
\end{center}\end{mydiagram}\end{minipage}
\hskip 1.5cm
\begin{minipage}{6cm}\begin{mydiagram}\label{fig8_quadruples}
\begin{center}
$$\includegraphics[width=6cm]{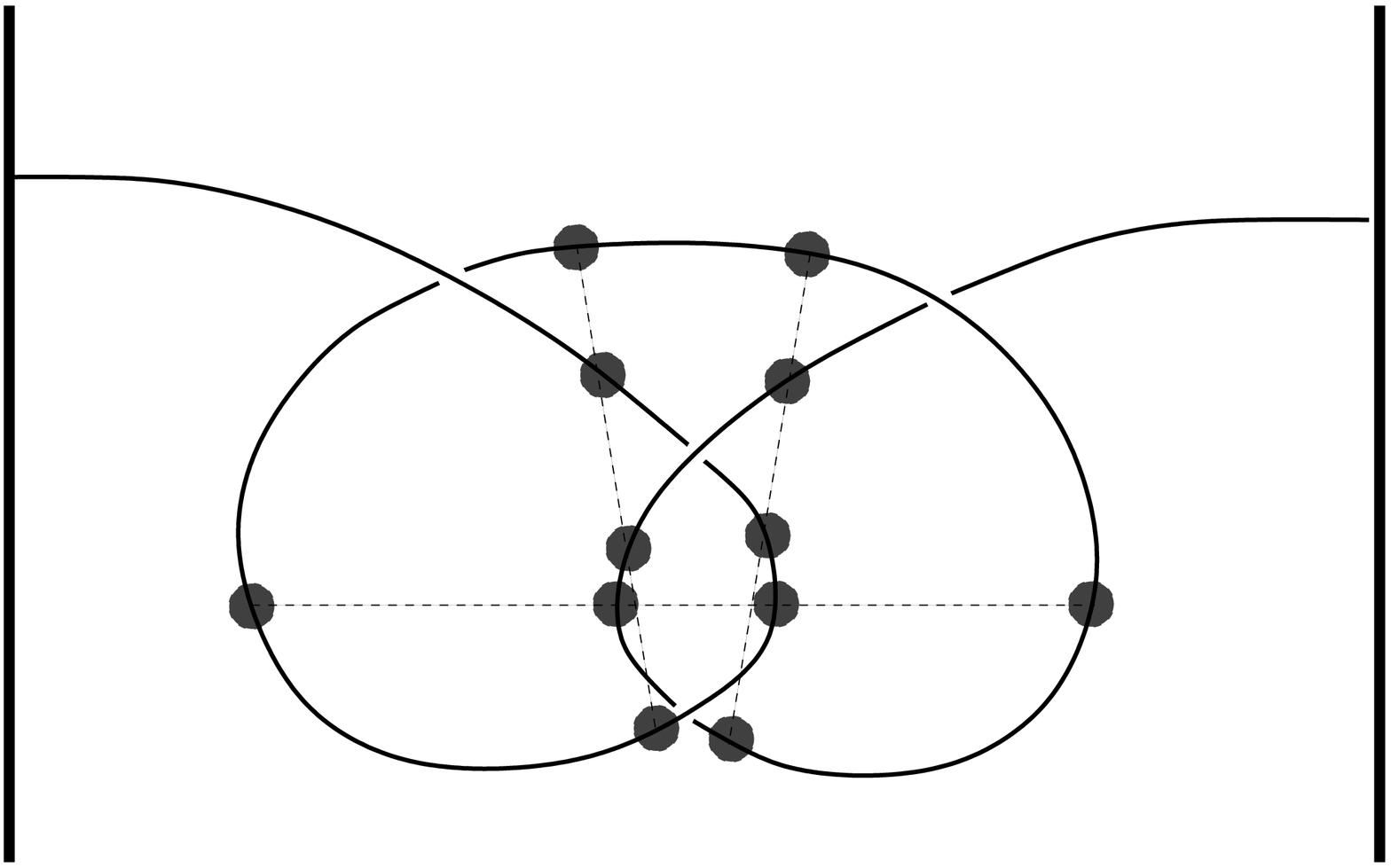}$$
quadrisecants.
\end{center}\end{mydiagram}\end{minipage}
\end{center}

Figures \ref{trefoil_quadruples} and \ref{fig8_quadruples} 
display all the quadrisecants on these knots.
For the trefoil, parametrize the knot so that the points where $x_3'(t)=0$
are $t \in \{ \frac{3}{12},\frac{5}{12},\frac{7}{12},\frac{9}{12} \}$.
These points partition $\I$ into five intervals.  Since $x_3$ is 
increasing on  $(\frac{3}{12},\frac{5}{12})$ and $(\frac{7}{12},\frac{9}{12})$,
and decreasing on $(0,\frac{3}{12})$, $(\frac{5}{12},\frac{7}{12})$ and 
$(\frac{9}{12},1)$, any quadrisecant must have its points on four 
of these five intervals, and no quadrisecant can have a point on the
$(\frac{5}{12},\frac{7}{12})$ interval for rather simple reasons -- such 
a quadrisecant of the projection could not lift to a quadrisecant of the knot itself
since the $x_3$ coordinate of such a  quadrisecant could not be monotone along the line.
Again by  monotonicity, there could be only one quadrisecant which passes through the
four intervals $[0,\frac{3}{12}]$,
$[\frac{3}{12},\frac{5}{12}]$,$[\frac{7}{12},\frac{9}{12}]$ and $[\frac{9}{12},1]$,
which is sketched in Figure \ref{trefoil_quadruples}.
%

\begin{center}
\begin{minipage}{7cm}\begin{mydiagram}\label{trefoil_simplex}
\begin{center}
{\psfrag{x}{$t_1$}\psfrag{y}{$t_2$}\psfrag{z}{$t_3$}
$$ \includegraphics[width=8cm]{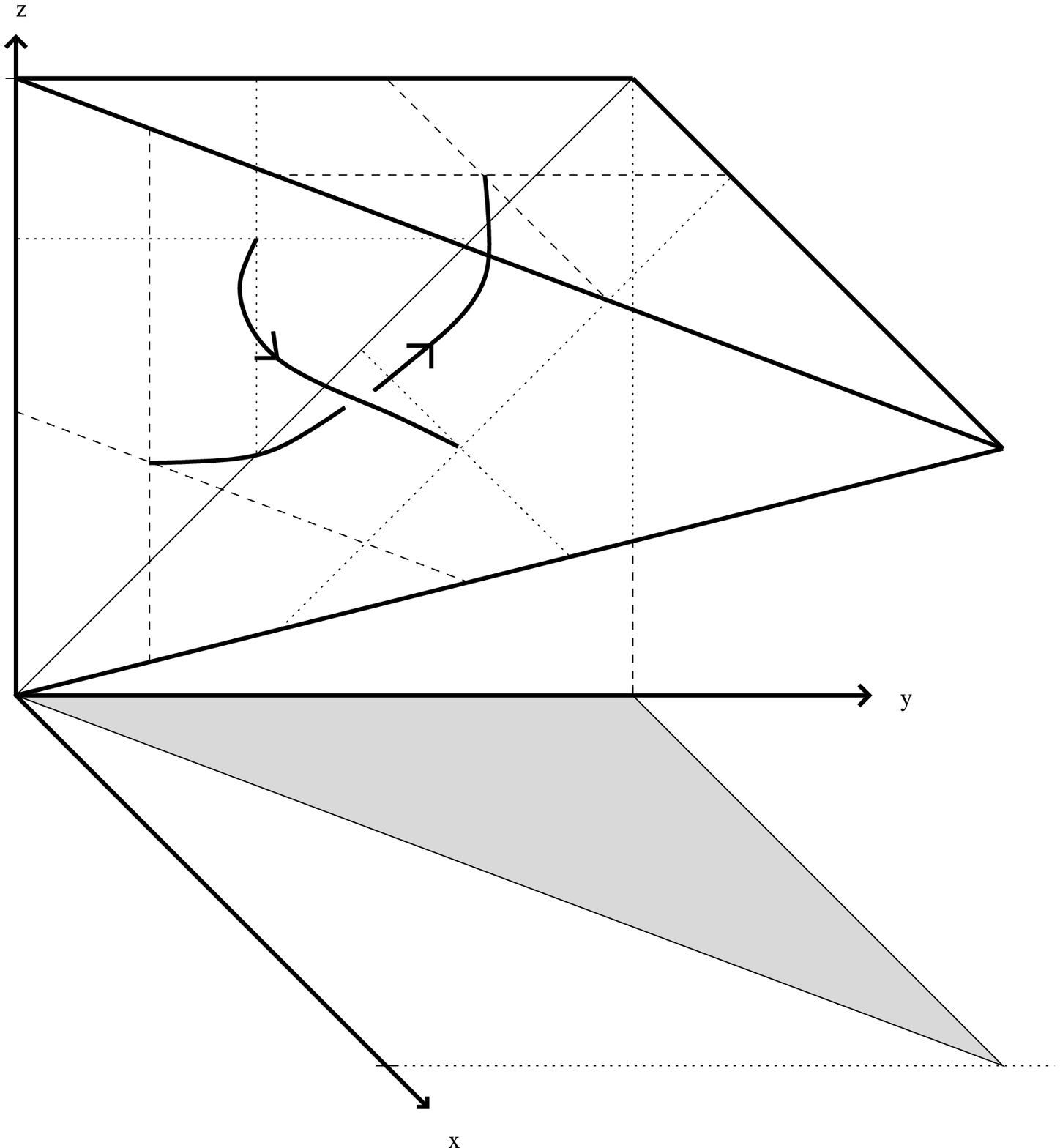}$$}
$\coll_1 [f] \cup \coll_3 [f] \subseteq \Delta^3$
\vskip 1mm
$f=\text{trefoil}$
\end{center}\end{mydiagram}\end{minipage}
\hskip 5mm
\begin{minipage}{7cm}\begin{mydiagram}\label{fig8_simplex}
\begin{center}
{\psfrag{x}{$t_1$}\psfrag{y}{$t_2$}\psfrag{z}{$t_3$}
$$ \includegraphics[width=8cm]{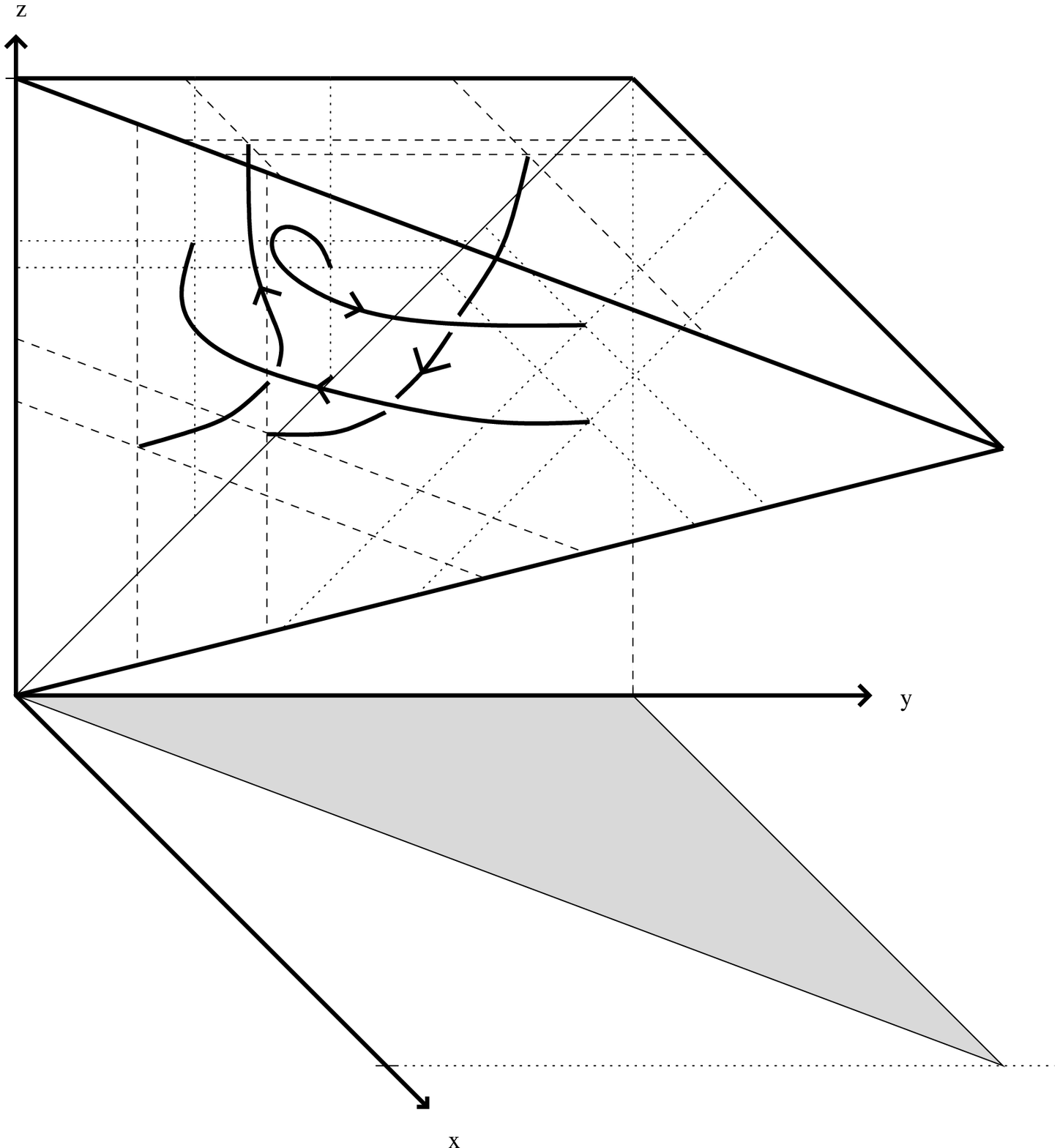}$$}
$\coll_1 [f] \cup \coll_3 [f] \subseteq \Delta^3$
\vskip 1mm
$f=\text{figure-8 knot}$
\end{center}\end{mydiagram}\end{minipage}
\end{center}

In Figures \ref{trefoil_simplex} and \ref{fig8_simplex} the data from 
the previous figures is assembled  to represent 
$\coll_i [f] \subseteq \Delta^3$ for $i=1,3$, together with orientation 
information which is explicitly computed below.  This
information suffices to compute our linking invariant. 
In the case of the trefoil, we know each of $\coll_i [f]$ for $i=1,3$ has 
only the one component, since we have
determined the boundary of these components and there are no closed 
components. We will not argue the non-existence of boundaryless 
components, but by Lemma~\ref{crossing} even if there were 
boundaryless components, they would be irrelevant since they would 
not cross any other components.
By following families of collinearities from the boundaries, one can see that 
none of the $\coll_i [f]$ has both boundary points on the same face of the 
simplex.  Moreover, in both the trefoil and figure eight examples, all the strand 
crossings in the diagrams Figures~\ref{trefoil_simplex} and \ref{fig8_simplex} 
are $\coll_3 [f]$ strands crossing over $\coll_1 [f]$ strands,
an elementary observation from Figures \ref{trefoil_quadruples} and 
\ref{fig8_quadruples}. 


To compute the orientations of the strands in Figures \ref{trefoil_simplex}
and \ref{fig8_simplex}, it suffices to compute the orientations induced on the
intersection of $\coll_i [f]$ with $\Delta^3_{(i=2)}$,
since every strand is incident to one of these faces.  
The orientation of these points is given 
(up to a choice of orientation of $\coll_i \I^3$) by the sign of
the determinant $\det[f(t_i)-f(t_k),f''(t_i),f'(t_k)]$
where $(t_1,t_2,t_3) \in \coll_i [f] \cap \Delta^3_j$ and $\{1,2,3\}=\{i,k,2\}$.
To see this, consider the identifications of the $(i=2)$ stratum of  
$C_3[\R^3]$ with $\R^3 \times (\R^3-\{0\}) \times S^2$. 
The orientation associated to $\coll_i[f]\cap \Delta^3_i$ is the transverse 
intersection number of $C_3[f]|_{\Delta^3_{(i=2)}}$ with 
$\coll_i[\I^3] \subseteq C_3[\I^3]$, which
quickly reduces to the sign of the above determinant.  

This analysis gives us the orientations in the figures, which allows us to 
establish the following.

\begin{theorem}\label{T:comp}
For the trefoil knot $\nu_2$  is $+1$, and for the figure-eight knot 
$\nu_2 = -1-1+1=-1$.
\end{theorem}

\section{Quadrisecants, finite-type, and geometric consequences} \label{coquad}

In the last section it was clear that the key to computing $\nu_2(\knot)$ was
identifying crossings of $\coll_1(\knot)$ and $\coll_3(\knot)$, which by
Lemma~\ref{crossing} correspond to quadrisecants on
the knot.  In this section we elaborate on this to compute $\nu_2(\knot)$
by counting quadrisecants of $\knot$, which leads to an easy proof that
$\nu_2$ coincides with the $z^2$-coefficient of the Conway polynomial.
We end this section with a result on
quadrisecants of a closed knot and applications to stick number of a knot.

\subsection{Quadrisecants on long knots}
The subspace of $C_4(\I^3)$ consisting of configurations of four points which
are \colinear\ has precisely twelve components.  If $(x_1,x_2,x_3,x_4)$ is a
quadrisecant, orient the line on which they sit by the convention
that the vector $x_2-x_1$ be positively oriented.  A choice of orientation
determines a permutation of $\{1,2,3,4\}$ given by $\sigma(i)=j$ if
the $i$-th point on the line is $x_j$. Note the permutations achieved are precisely 
the permutations such that $\sigma(2)>\sigma(1)$, and there are 
twelve of these.

\begin{definition}\label{signdef}
Let $\CO$ denote the subspace of $C_4 (\R^3)$  of \colinear\
configurations labeled by the $4$-cycle $(1 3 4 2)$.
Let $\knot \subseteq \I^3$ be a knot, parametrized by $f : \I \to \I^3$.
We associate a sign $\epsilon_x$ to a quadruple $x=(f(t_1),f(t_2),f(t_3),f(t_4)) \in
\CO$ by defining it to be the sign of the determinant of the $2 \times 2$ matrix:
$$
\left[
\begin{matrix}
|f(t_3) - f(t_2)|\cdot \det[v, f'(t_1),f'(t_3)] & |f(t_3)-f(t_1)| \cdot \det[v, f'(t_2),f'(t_3)] \\
|f(t_4) - f(t_2)|\cdot \det[v, f'(t_4),f'(t_1)] & |f(t_4)-f(t_1)| \cdot \det[v, f'(t_2),f'(t_4)] \\
\end{matrix}
\right]
$$
\noindent where $v=f(t_2)-f(t_1)$.
\end{definition}

With this sign convention in hand, we give an alternate definition of our 
self-linking invariant.

\begin{proposition}\label{co4prop}
Let $\knot = \im(f)$ be a knot $\knot \subseteq \I^3$.
If  $C_4[\knot]$ is transverse to $\CO$, then
$$ \nu_2(\knot) = \sum_{x \in C_4(\knot) \cap \CO} \epsilon_x$$
\end{proposition}

\begin{proof}
As in Lemma~\ref{crossing}, project the link 
$C_3[f]^{-1}(\coll_i [\I^3])$ onto $\Delta^2$ by $\rho$, which forgets 
the $t_1$ coordinate.  Let $L = \overline{\coll_1[f]} \cup \overline{\coll_3[f]}$.
Because $\coll_1 [f]$ has boundary in $\Delta^3_{(1=2)}$ and 
$\Delta^3_{(3=4)}$, and $\coll_3 [f]$ has boundary in 
$\Delta^3_{(2=3)}$ and $\Delta^3_{(1=0)}$,  the crossing 
points of the $\rho$-projection of  $L$ differ from the crossing points of the 
$\rho$-projection of $\coll_1[f] \cup \coll_3[f]$ by adding crossings where 
$\coll_1[f]$ strands cross over $\coll_3[f]$ strands.  Therefore the
linking number of $L$ is precisely the crossing number of the $\rho$-projection
of $\coll_1[f] \cup \coll_3[f]$, where we count only
$\coll_3[f]$ strands crossing over $\coll_1[f]$.

We next investigate such crossings.  As in Lemma~\ref{crossing},
denote the two preimages of the crossing point by $(t_1^*, t_2^*, t_3^*)$
and $(t_1', t_2' = t_2^*, t_3' = t_3^*)$.  We deduce, by the definitions
of $\coll_i[f]$ and the fact that $\coll_3[f]$ is the over strand, that
 $t^*_1 > t'_1$ with $f(t^*_3)$ between $f(t^*_1)$ and $f(t^*_2)$ and
 $f(t'_1)$ between $f(t^*_2)$ and $f(t^*_3)$ on the same line $L$ .  With
 respect to $L$'s orientation we have the order relation
$f(t^*_1)<f(t^*_3)<f(t^*_2)$ so the only possible order relation for the
$f(t'_i)$'s on $L$ is $f(t'_3)<f(t'_1)<f(t'_2)$ since only the $t_1$
coordinate changes. Therefore, we have the order relation
$f(t^*_1)<f(t^*_3)<f(t'_1)<f(t^*_2)$ on $L$ and thus our permutation is
$\sigma_x = (1 3 4 2)$.

The sign of $\epsilon_x$ is straightforward to justify.
As argued above, we are counting the transverse intersection number
$C_4[\knot] \cap \CO$. Orient the image of $C_4[\knot]$ using the standard product
orientation of its domain and  $\CO$ by making the identification
$\R^3 \times (\R^3 -\{0\}) \times (0,\infty) \times (0,1) \equiv \CO$
by $(x,v,t_1,t_2) \longmapsto (x,x+v,x-t_1v,x+t_2v)$
and giving it the product orientation.   
Both the image of $C_4[\knot]$ and $\CO$ have trivial tangent bundles, given by
the above product structures. Thus, the points in the transverse
intersection $C_4[\knot]\cap \CO$ have orientations given by the sign of
the determinant of the basis consisting of the two trivializations,
$$ \det\left\{ \begin{matrix}
I_{\R^3} & 0        & 0        & 0        & f'(t_1) & 0       & 0 & 0 \\
I_{\R^3} & I_{\R^3} & 0        & 0        & 0       & f'(t_2) & 0 & 0 \\
I_{\R^3} & -\frac{|f(t_3)-f(t_1)|}{|f(t_2)-f(t_1)|}I_{\R^3} & f(t_1)-f(t_2) & 0        & 0       & 0       & f'(t_3) & 0 \\
I_{\R^3} & \frac{|f(t_4)-f(t_1)|}{|f(t_2)-f(t_1)|}I_{\R^3} & 0 & f(t_2)-f(t_1) & 0       & 0       & 0       & f'(t_4)
\end{matrix}
\right\}$$
which rapidly reduces to the formula in the statement of the proposition.
\end{proof}

\subsection{$\nu_2$ is of type two.}

\begin{theorem}\label{T:typetwo}
$\nu_2$ is of finite type two.
\end{theorem}
\begin{proof}
We show the invariant is type two by computing directly that its
third derivative (in the sense of finite-type invariant theory) is zero.
We show that for every knot $\knot$ and a set of three crossing
changes $c_1,c_2,c_3$,
that the associated alternating sum vanishes:
\begin{align}
\sum_{\sigma\subset [3]}(-1)^{|\sigma |} \nu_2(\knot_\sigma)=0, \label{altsum}
\end{align}
where $[3]$ denotes the set $\{1,2,3\}$ and $\knot_\sigma$ denotes the knot
obtained from $\knot$ by applying crossing changes $c_i$ where $i \in \sigma$.
The three crossing changes are
supported in balls $B_i, 1\leq i\leq 3$.
It suffices to check this when the three balls are not collinear, i.e. no straight line
intersects all three balls, since an isotopy can always bring us to that situation.

Given a quadrisecant on one of the eight knots $\knot_\sigma$, we know that it avoids one
of the balls, say $B_3$. Thus a quadrisecant appears in the knot $\knot_\tau$ where
$\tau\subset [2]$ if and only if it appears in $\knot_{\tau\cup\{3\}}$. Taking the signed
sum, the net contribution is zero.

\end{proof}


\begin{corollary}\label{C:n2=c2}
The self-linking invariant, $\nu_2 (\knot)$ is equal to the $z^2$ coefficient of the
Alexander-Conway polynomial, $c_2 (\knot)$.
\end{corollary}

\begin{proof}
The space of type two invariants is spanned by the
Conway coefficient $c_2$ and the constant function $1$, which have
independent values on the unknot and the trefoil. 
By \refT{typetwo}, $\nu_2$ is type two.  Both $\nu_2$ and $c_2$
vanish on the unknot  and take value $1$ on the trefoil by the previous section, and
therefore must coincide.
\end{proof}

\begin{remark}
It is straightforward from definitions that counting
quadrisecants on the embedding level is type $2$, which led Bar-Natan,
Hutchings, and D. Thurston to conjecture a corresponding invariant 
at the level of isotopy.
\end{remark}

\subsection{Quadrisecants and closed knots.}

\begin{definition}
Let $\overline{\knot}$ denote the convex hull of the knot and
define the set of extremal points of the knot to be
$\knot \cap \partial \overline{\knot}$, generically a finite number of closed 
intervals provided the knot is non-trivial..
\end{definition}

Given an oriented closed knot $\knot$ with parametrization 
$f : S^1 \to \R^3$ and an extremal point $p \in \knot$, 
linearly order the points in $\knot - \{p\}$ 
using the ordering induced by the orientation of $\knot$.
Consider the component of $C_4[\knot - \{p\}]$ for which
$t_1 \leq t_2 \leq t_3 \leq t_4$, which we call $C^o_4[\knot - \{p\}]$.  
Provided $C_4[f] : C^o_4[S^1 - f^{-1}(p)] \to C_4[\R^3]$
is transverse to $\CO$, we can give a point in the intersection 
$x \in C_4(\knot) \cap \CO$  such that $x=(f(t_1),f(t_2),f(t_3),f(t_4))$ a 
sign $\epsilon_x$ exactly as in  Definition~\ref{signdef}.  

\begin{lemma}\label{roundlong}
With conventions as above,
$$ c_2 \knot = \sum_{x \in C_4(\knot) \cap \CO} \epsilon_x.$$
\end{lemma}

\begin{proof}
Translate and scale the knot so that its image lies in  $\I^3$. 
``Open'' the knot at $p$ to give a long knot $\tilde f : \I \to \I^3$.
 When we open the knot to create $\tilde f$ we may create new
quadrisecants, but by design the associated permutation $\sigma_x$ of
these quadruples must fix either $1$ or $4$.
Thus the sum $\sum_{x \in C_4(\knot) \cap \CO} \epsilon_x$ is
is precisely $\nu_2 (\tilde f)$, which by \refC{n2=c2} is $c_2(\tilde f)$.
\end{proof}

We now establish that the $(1342)$-quadrisecants of knots with boundary 
correspond to a special kind quadrisecant of closed knot.

\begin{definition}\label{D:nsns}
Let $L$ be an oriented line that intersects
an oriented knot $\knot$ in four points $P= \{p_0, \ldots, p_3\} = \{x_1,\ldots,x_4\}$,
where the subscripts of $x_i$'s are determined by the orientation of $L$ and
the subscripts of $p_i$'s (which are understood modulo four)
are consistent with the cyclic ordering given by $\knot$.  Let $m = \frac{x_2 + x_3}{2}$.
We call $P$ an {\em alternating quadrisecant} 
if $p_i - m$ and $p_{i+1} - m$ are  
negative multiples of each other for all $i$.
\end{definition}

 This definition is invariant under the action  the subgroup
$\langle (1234),(14)(23)\rangle \subset S_4$ on indices.

\begin{definition}
Let $\knot$ and $L$ intersect in 
an alternating quadrisecant $P$ as above. We call the component of
$\knot \setminus P$ that has boundary equal to $\{x_2, x_3\}$
the middle component of the quadrisecant on $\knot$.
\end{definition}

\begin{center}
\begin{minipage}{9cm}\begin{mydiagram}\label{externalarc}
\begin{center}
{
\psfrag{C}{$C$}
\psfrag{l1}{$p_0 = x_3$}
\psfrag{l2}{$p_1 = x_1$}
\psfrag{l3}{$p_2 = x_4$}
\psfrag{l4}{$p_3 = x_2$}
\psfrag{L}{}
$$\includegraphics[width=9cm]{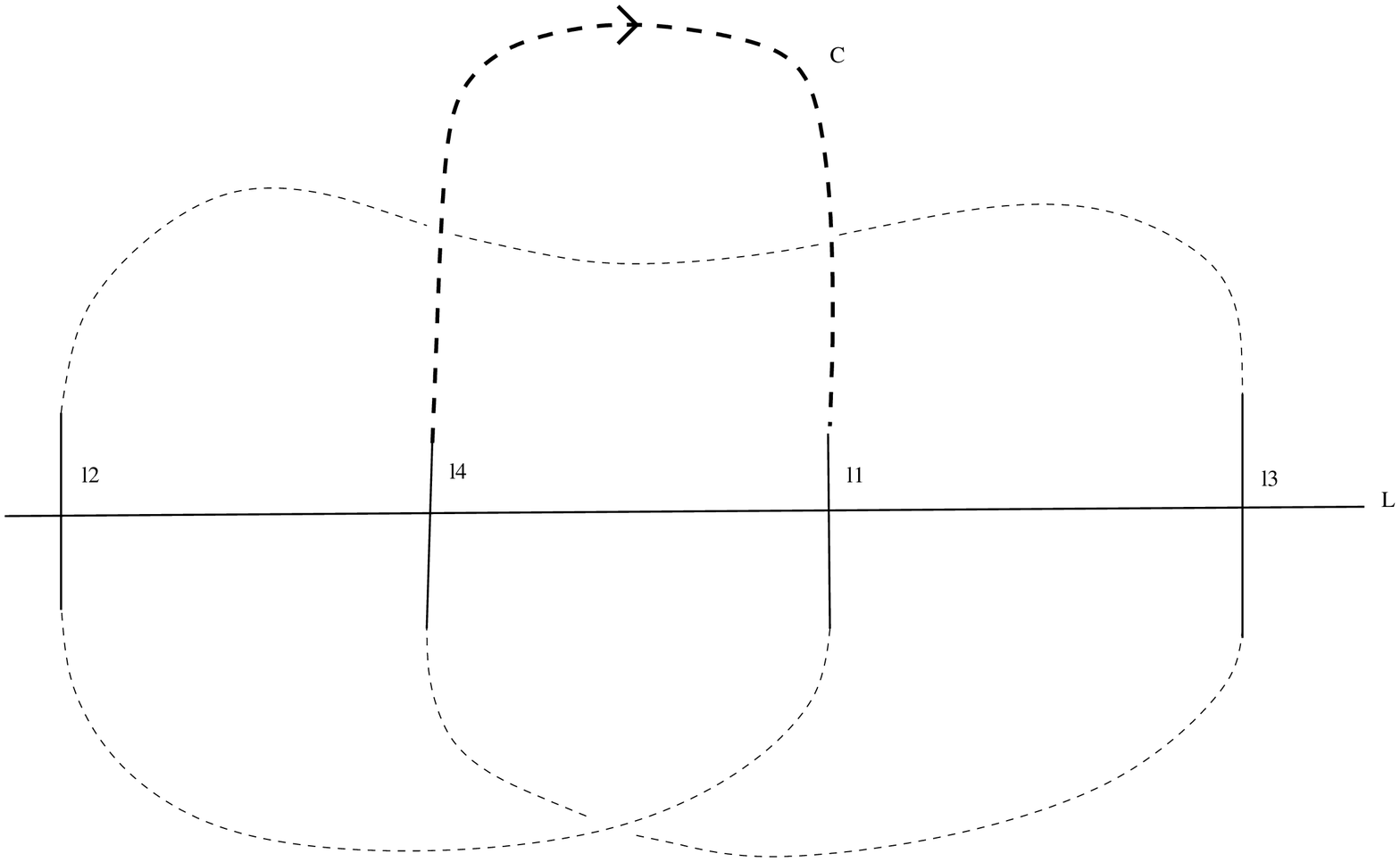}$$}
A quadrisecant $P = \knot \cap L$ with middle component $C$.
\end{center}\end{mydiagram}\end{minipage}
\end{center}

\begin{lemma}\label{nsnsrec}
Given a quadrisecant $P = L \cap \knot$ of a closed knot $\knot$, it is a
$(1342)$-labeled quadrisecant for some extremal point $p \in \knot$ if and only if the quadrisecant is alternating and its middle component $C$ contains
an extremal point of the knot.
\end{lemma}

\begin{proof}
If $p$ is an extremal point making $x$ into a $(1342)$-labeled quadrisecant,
it is clear that the quadrisecant is alternating.
If, on the other hand, $L$ is a line intersecting $\knot$ in an alternating quadrisecant,
and $p \in C$ is an extremal point, begin assigning a cyclic ordering to the
points of $x$ by giving the label $1$ to
the first point of $x$ that occurs after $p$ in the ordering of $\knot$.  The
label $4$ must be between the points labeled $1$ and $2$ as this is a cyclic ordering and
$2$ and $4$ must occur on the same side of the midpoint $m$.  Thus the permutation
associated to this quadrisecant is $\sigma_x = (1342)$.
\end{proof}

The configuration space $C_4[S^1]$ has six components; let $C_4^o[S^1]$ denote the
component where a configuration $(t_1,t_2,t_3,t_4)$ has the cyclic ordering
$t_1 \leq t_2 \leq t_3 \leq t_4 \leq t_1$.  Let us also be more definite about
the labeling of points in an alternating quadrisecant $P = L \cap \knot$ from \refD{nsns}.
Choose $\{ p_0, \ldots, p_3\} = P$ so that the boundary of the middle component  $C$
of $P$ on $\knot$ is $\{ +p_0, -p_3\}$, as in Figure~\ref{externalarc}.
Let $v = p_1 - p_0$, and let $p'_i$ be the unit tangent vector of $\knot$ at $p_i$.
Also, we may assume $\knot \cap \partial \overline{\knot}$ has a
finite number of components.  We let
$n_L$ be the number of components of $C \cap \partial \overline{\knot}$.

\begin{proposition}\label{co4closed}
Given an oriented knot $\knot=\im(f)$, such that 
$C_4[f] : C_4^o[S^1] \to C_4[\R^3]$
is transverse to $\CO$, and $\knot \cap \partial \overline{\knot}$ is a 
finite collection of intervals then
$$c_2 \knot = \frac{1}{|\pi_0 (\knot \cap \partial \overline{\knot}) |}
\cdot \sum_L n_L sig \left(\det
 \left[
\begin{matrix}
|p_2 - p_1|\cdot \det[v, p'_0,p'_2] & |p_2-p_0| \cdot \det[v, p'_1,p'_2] \\
|p_3 - p_1|\cdot \det[v, p'_3,p'_0] & |p_3-p_0| \cdot \det[v, p'_1,p'_3] \\
\end{matrix}
\right]\right)
 $$
where the sum is taken over lines $L$ that intersect $\knot$ in an alternating quadrisecant.
\end{proposition}

\begin{proof}
Combining Lemmas \ref{roundlong} and \ref{nsnsrec}, the above sum is simply
$\frac{1}{|\pi_0 (\knot \cap \partial \overline{\knot} )|}$ times the sum
of the equation given in Lemma \ref{roundlong} where we sum over the
components of $\pi_0 \knot \cap \partial \overline{\knot}$, choosing one
extremal point for each component of $\knot \cap \partial
\overline{\knot}$.
\end{proof}

\begin{corollary} \label{co4cor}
Given a knot $S^1 \to \R^3$ that has a non-zero $z^2$ coefficient in its
Alexander-Conway polynomial, then there exists an alternating quadrisecant
of the knot
\end{corollary}

Proposition~\ref{co4closed} gives a partial converse to the main
results of \cite{Pa33, MM82, Ku94}, in which a knot 
is shown to be trivial if it has no quadrisecants.  
Corollary \ref{co4cor} resolves a special case of the conjecture of 
Cantarella, Kuperberg, Kusner, and Sullivan that alternating quadrisecants 
exist for all non-trivial knots \cite{CKKS02}.
The full conjecture has recently been resolved 
in the dissertation of Denne \cite{De03}, which is currently being written.

\subsection{Geometric applications}

We give simple applications of these results to bounding the complexity
of minimal polygonal and
polynomial realizations of a knot. The stick number of a knot is the minimum number
of line segments among all piecewise linear knots equivalent to a given one in $\R^3$. 
There are various relationships known between the stick number of a knot and
other invariants of the knot which are difficult to compute such as the crossing number.
See \cite{RS03} for a survey of known results on the stick number of knots.

Let $K$ be a polygonal knot.
If $L$ is a straight line that intersects $K$ in precisely four points, then
all the points of $K \cap L$ sit on different line segments in $K$. 
We associate to $L$ a set $seg(L)$ consisting of the four line segments of 
$K$ that contain $K \cap L$.

\begin{lemma}\label{countinglemma}
A polygonal knot in general position has finitely many quadrisecants.  
The function $L \longmapsto seg(L)$, is at most two-to-one.
\end{lemma}

\begin{proof}
Let $f_i(t)=x_i+tv_i$ for some $x_i \in R^3$ and $v_i \in \R^3-\{0\}$ and $t \in \I$
be parametrizations of  line segments of $K$ for $i \in \{1,2,3,4\}$.
Apply a small isotopy of $K$ so that any three of the four vectors
$v_1,v_2,v_3,v_4$ are linearly independent.  Now consider the extensions of
$f_i$ to all of $\R$, and restrict our attention to the first three
functions $f_1,f_2,f_3$. There are three cases to consider, as 
the images of the $f_i$ are disjoint, or one pair, say $f_1$ and $f_2$ intersect in a point, 
or  $f_i$ intersects $f_{i+1}$ for $i \in \{1,2\}$.

In the first case there is a unique affine-linear transformation $T$ of $\R^3$ such
that $T \circ f_1(t)=(t,0,0)$, $T \circ f_2(t)=(0,t,1)$ and $T \circ f_3(t)=(1,1,t)$.
Affine-linear transformations do not change collinearity properties, so without
loss of generality we can restrict ourselves to studying  quadrisecants
of these three lines and an arbitrary $f_4(t)=x_4+tv_4$. 
The closure of the union of all lines $L$ which intersect each of $\im(f_i)$ 
for $i\in \{1,2,3\}$
is the algebraic variety $V=\{(x,y,z)\in \R^3 : (x+z-1)y=xz\}$. 
Moreover,  every point in $V$  belongs to at most one
line that intersects  each of $\im(f_i)$ for $i\in \{1,2,3\}$. Therefore, 
a straight line $L$ that intersects $\im(f_i)$ for $i \in \{1,2,3,4\}$ corresponds
to a point in the intersection $V \cap \im(f_4)$.  Because
$V$ is a quadratic surface  $V \cap \im(f_4)$ can consist of
$0,1,2$ or an infinite number of points.  
Since $V$ is a ruled surface, infinite intersections
$V \cap \im(f_4)$ occur only when $\im(f_1) \subset V$ which is  a codimension
three condition, not occurring generically. This proves the theorem in case 1. 

The last two cases are almost identical to the first.  In the second, we can use an affine-linear
transformation to convert the straight lines $f_1,f_2,f_3$ into the family
$f_1(t)=(t,0,0)$, $f_2(t)=(0,t,0)$ and $f_3(t)=(1,1,t)$.  In this case the algebraic
variety of collinear triples is a subvariety of $\{(x,y,z) : x=y \text{ or } z=0 \}$, and
the result follows immediately.  In the third case, we use an affine linear transformation
to convert the first three line segments to $f_1(t)=(t,0,0)$, $f_2(t)=(0,t,0)$ and
$f_3(t)=(0,1,t)$, making the variety of collinear 
triples $\{(x,y,z) : z=0 \text{ or } x=0 \}$.
\end{proof}

\begin{theorem}
Given a PL knot K, let $n(K)$ be the number of line segments in $K$. Then 
$$\frac{|c_2|}{2} \leq {n(K)\choose 4}$$
\end{theorem}

\begin{proof}

If $K$ is in general position, we can assume it has a finite number 
of quadrisecants, and that no quadrisecant has a point in the
$0$-skeleton $K^0$ of $K$. 
Let $K_\epsilon$ be the epsilon-smoothing of the knot $K$, which is 
a smooth knot such that if $B_{\epsilon}(v)$ is the epsilon ball about 
a vertex $v$ of $K$, then 
$K-\cup_{v \in K^0} B_\epsilon(v)  = K_\epsilon - \cup_{v \in K^0} B_\epsilon(v)$
and the knots $K \cap B_\epsilon(v)$ and $K_\epsilon \cap B_\epsilon(v)$
are isotopic in a restricted sense. 


For some $\epsilon >0$, $K_\epsilon$ has all the quadrisecants
of $K$, and we leave it to the reader
to argue that it can have no more alternating quadrisecants
using compactness of the segments of $K$ and the fact that none of the
quadrisecants intersect $K^0$. 
Because $K_{\epsilon}$ satisfies the criteria of Lemma \ref{roundlong},
$K$ satisfies the
criteria of Lemma \ref{countinglemma}, and both knots have exactly the same collection
of alternating quadrisecants,  the inequality
$\frac{|c_2|}{2} \leq {n(K)\choose 4}$
comes immediately from the fact that the function $L \longmapsto seg(L)$
is at most two to one.
\end{proof}

In a similar vein is the following theorem.

\begin{theorem}
Let $\knot$ be a long knot $\R^1 \to \R^3$ parametrized by polynomials 
of degree $n$ with leading coefficient $1$.  If 
$C_4[\knot]$ is transverse to $\CO$ then $|c_2(\knot)| \leq (2n)^4$.
\end{theorem}

\begin{proof}
If $x(t)$, $y(t)$ and $z(t)$ parametrize the knot, then a collinearity at times
$t_1, \ldots, t_4$ translates to a solution to a system of equations including
for example
$$\left(x(t_1) - x(t_2)\right)\left(y(t_2) - y(t_ 3)\right) = 
\left(y(t_1) - y(t_2)\right)\left(x(t_2) - x(t_ 3)\right).$$
There are four such equations, each of degree $2n$ so Bezout's theorem
implies that if there are finitely many solutions there are at most $(2n)^4$.
\end{proof}

Using Propositions ~\ref{co4prop} and \ref{co4closed} to find relationships
between $|c_2|$ and more subtle and geometric but less computable invariants
is worthy of further study.  The bounds we have just obtained for polygonal
and polynomials knots are likely not the best possible.
It also seems reasonable that one might be able to get a lower-bound on 
the total curvature  of a knot using the above results on quadrisecants for a knot.  
According to Istv\'an F\'ary \cite{Fa49}, perhaps the first proof of the 
F\'ary-Milnor theorem is due to Heinz Hopf. 
F\'ary mentions that the key part of Hopf's unpublished
proof is the quadrisecant result of Pannwitz.

We conclude our paper with illustrations of Proposition \ref{co4closed}.

\begin{center}
\scalebox{0.6}{\includegraphics{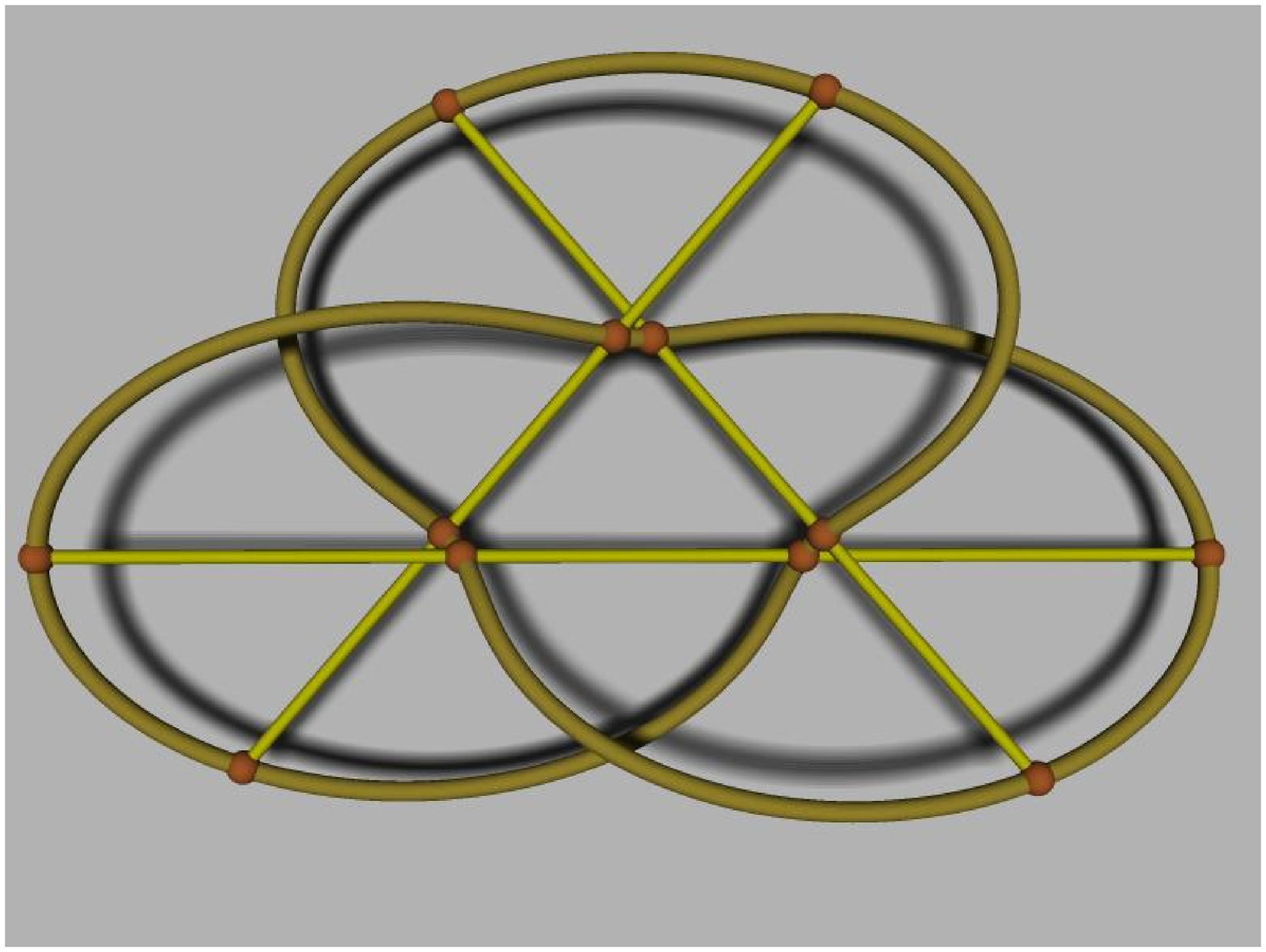}}

All $\epsilon_L=+1$, $n_L=1$, $\pi_0 \knot \cap \partial \overline{\knot} = 3$
\end{center}

\begin{center}
\scalebox{0.6}{\includegraphics{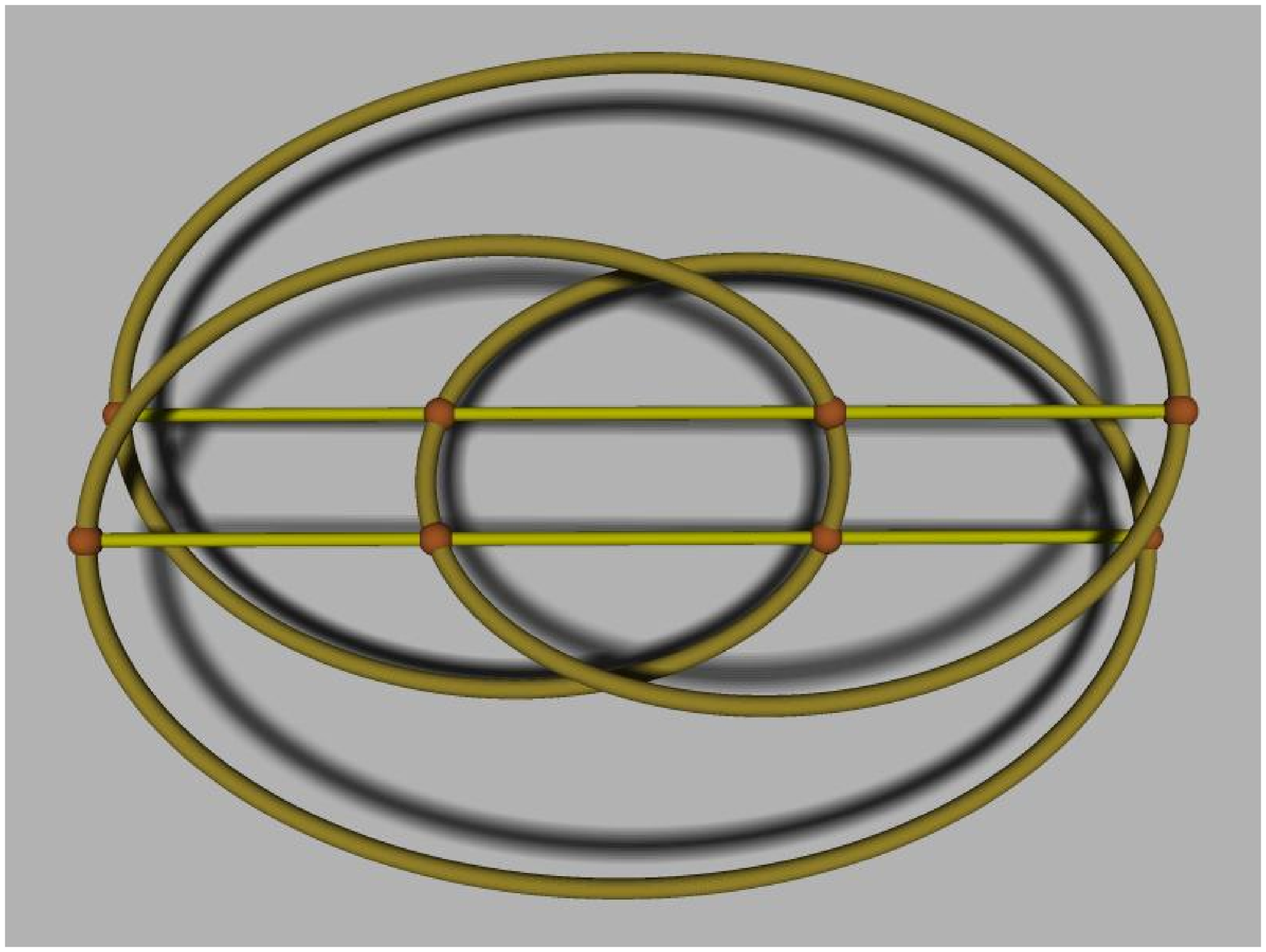}}

All $\epsilon_L=+1$, $n_L=1$, $\pi_0 \knot \cap \partial \overline{\knot} = 2$
\end{center}

\begin{center}
\scalebox{0.6}{\includegraphics{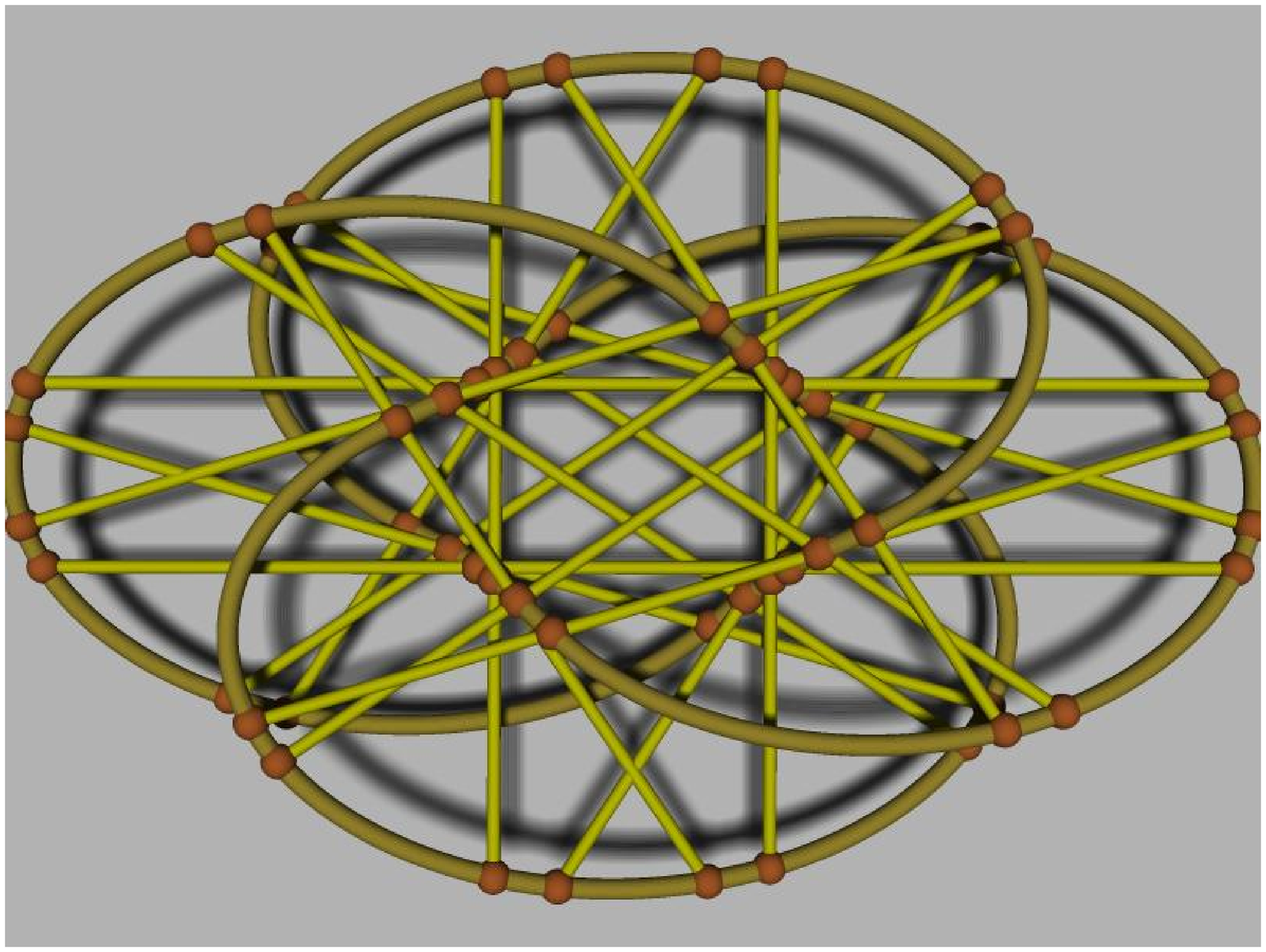}}

All $\epsilon_L=+1$, $4$ $L$ have $n_L=2$ and $12$ $L$ have $n_L=1$, $\pi_0 \knot 
\cap \partial \overline{\knot} = 4$
\end{center}

\begin{center}
\scalebox{0.6}{\includegraphics{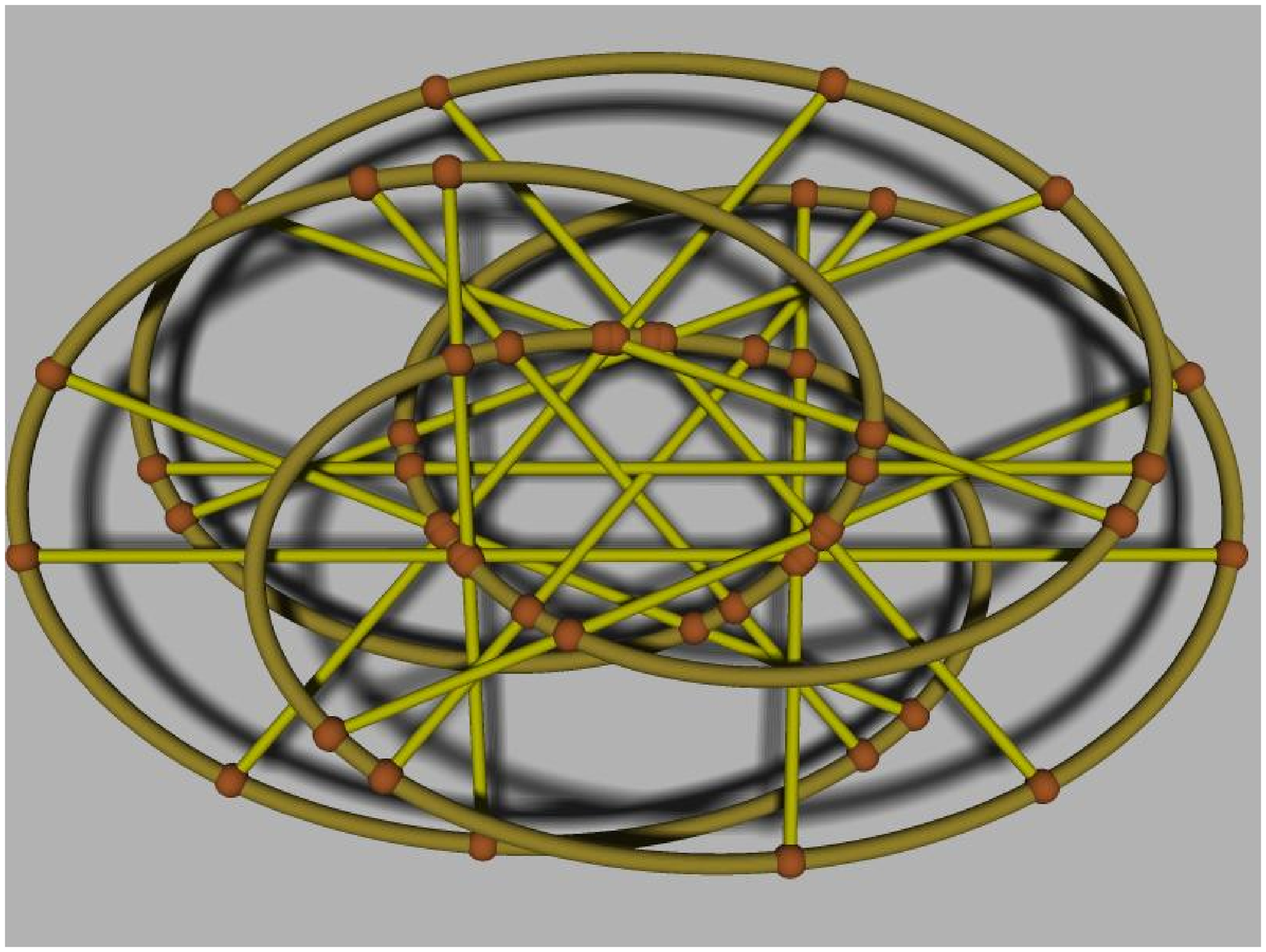}}

All $\epsilon_L=+1$, $3$ $L$ have $n_L=2$ and $9$ $L$ have $n_L=1$, 
$\pi_0 \knot \cap \partial \overline{\knot} = 3$
\end{center}

\bibliographystyle{amsplain}
\bibliography{data}

\providecommand{\bysame}{\leavevmode\hbox to3em{\hrulefill}\thinspace}
\providecommand{\MR}{\relax\ifhmode\unskip\space\fi MR }
\providecommand{\MRhref}[2]{%
  \href{http://www.ams.org/mathscinet-getitem?mr=#1}{#2}
}
\providecommand{\href}[2]{#2}
\begin{thebibliography}{10}

\bibitem{AF97}
D.~Altsch\"uler and L.~Freidel, \emph{Vassiliev knot invariants and
  {C}hern-{S}imons perturbation theory to all orders}, Comm. Math. Phys.
  \textbf{187} (1997), no.~2, 261--287.

\bibitem{AS94}
S.~Axelrod and I.~M. Singer, \emph{Chern-{S}imons perturbation theory. {I}{I}},
  J. Differential Geom. \textbf{39} (1994), no.~1, 173--213.

\bibitem{BN95}
D.~Bar-Natan, \emph{On the {V}assiliev knot invariants}, Topology \textbf{34}
  (1995), no.~2, 423--472.

\bibitem{Bi93}
J.~S. Birman, \emph{New points of view in knot theory}, Bull. Amer. Math. Soc.
  (N.S.) \textbf{28} (1993), no.~2, 253--287.

\bibitem{BL93}
J.~S. Birman and X.-S. Lin, \emph{Knot polynomials and {V}assiliev's
  invariants}, Invent. Math. \textbf{111} (1993), no.~2, 225--270.

\bibitem{BT94}
R.~H. Bott and C.~H. Taubes, \emph{On the self-linking of knots}, J. Math.
  Phys. \textbf{35} (1994), no.~10, 5247--5287.

\bibitem{BCSS03}
R.~D. Budney, J.~Conant, K.~P. Scannell, and D.~P. Sinha, \emph{Mapping space
  models and finite-type knot invariants}, In preparation.

\bibitem{CKKS02}
J.~Cantarella, G.~Kuperberg, R.~B. Kusner, and J.~M. Sullivan, \emph{The second
  hull of a knotted curve}, To appear, Amer. J. Math., 2002.

\bibitem{CT03}
J.~Conant and P.~Teichner, \emph{Grope cobordism of classical knots}, to appear
  in Topology, math.GT/0012118, 2000.

\bibitem{De03}
E.~Denne, Ph.D. thesis, University of Illinois, Urbana-Champaign, 2003.

\bibitem{DK83}
W.~G. Dwyer and D.~M. Kan, \emph{Function complexes for diagrams of simplicial
  sets}, Nederl. Akad. Wetensch. Indag. Math. \textbf{45} (1983), no.~2,
  139--147.

\bibitem{FH01}
E.~R. Fadell and S.~Y. Husseini, \emph{Geometry and topology of configuration
  spaces}, Springer-Verlag, New York-Berlin-Heidelberg, 2001.

\bibitem{Fa49}
I.~F\'ary, \emph{Sur la courbure totale d'une courbe gauche faisant un noeud},
  Bull. Soc. Math. France \textbf{77} (1949), 128--138.

\bibitem{GK00}
T.~G. Goodwillie and J.~R. Klein, \emph{Excision statements for spaces of
  embeddings}, In preparation.

\bibitem{GKW01}
T.~G. Goodwillie, J.~R. Klein, and M.~S. Weiss, \emph{Spaces of smooth
  embeddings, disjunction and surgery}, Surveys on surgery theory, vol. 2
  (S.~E. Cappell, A.~A. Ranicki, and J.~M. Rosenberg, eds.), Ann. of Math.
  Stud., vol. 149, Princeton Univ. Press, Princeton, 2001, pp.~221--284.

\bibitem{GW99a}
T.~G. Goodwillie and M.~S. Weiss, \emph{Embeddings from the point of view of
  immersion theory. {I}{I}}, Geom. Topol. \textbf{3} (1999), 103--118.

\bibitem{Gou94}
M.~N. Goussarov, \emph{On $n$-equivalence of knots and invariants of finite
  degree}, Topology of manifolds and varieties (O.~Ya. Viro, ed.), Adv. Soviet
  Math., vol.~18, Amer. Math. Soc., Providence, 1994, pp.~173--192.

\bibitem{GP74}
V.~Guillemin and A.~Pollack, \emph{Differential topology}, Prentice-Hall, Inc.,
  Englewood Cliffs, 1974.

\bibitem{Ha00}
K.~Habiro, \emph{Claspers and finite type invariants of links}, Geom. Topol.
  \textbf{4} (2000), 1--83.

\bibitem{Hi55}
P.~J. Hilton, \emph{On the homotopy groups of the union of spheres}, J. London
  Math. Soc. (2) \textbf{30} (1955), 154--172.

\bibitem{Hi03}
P.~S. Hirschhorn, \emph{Model categories and their localizations}, Math.
  Surveys Monographs, vol.~99, Amer. Math. Soc., Providence, 2003.

\bibitem{Ko93}
M.~L. Kontsevich, \emph{Vassiliev's knot invariants}, I. {M}. {G}elfand
  {S}eminar (S.~Gelfand and S.~Gindikin, eds.), Adv. Soviet Math., vol.~16,
  Amer. Math. Soc., Providence, 1993, pp.~137--150.

\bibitem{Ko99}
M.~L. Kontsevich, \emph{Operads and motives in deformation quantization}, Lett.
  Math. Phys. \textbf{48} (1999), no.~1, 35--72.

\bibitem{Ko97}
U.~Koschorke, \emph{A generalization of {M}ilnor's $\mu$-invariants to
  higher-dimensional link maps}, Topology \textbf{36} (1997), no.~2, 301--324.

\bibitem{Ku94}
G.~Kuperberg, \emph{Quadrisecants of knots and links}, J. Knot Theory
  Ramifications \textbf{3} (1994), no.~1, 41--50.

\bibitem{KT99}
G.~Kuperberg and D.~P. Thurston, \emph{Perturbative $3$-manifold invariants by
  cut-and-paste topology}, math.GT/9912167, 1999.

\bibitem{MP77}
R.~S. Millman and G.~D. Parker, \emph{Elements of differential geometry},
  Prentice-Hall, Inc., Englewood Cliffs, 1977.

\bibitem{MM82}
H.~R. Morton and D.~M.~Q. Mond, \emph{Closed curves with no quadrisecants},
  Topology \textbf{21} (1982), no.~3, 235--243.

\bibitem{Pa33}
E.~Pannwitz, \emph{Eine elementargeometrische {E}igenshaft von
  {V}erschlingungen und {K}noten}, Math. Ann. \textbf{108} (1933), 629--672.

\bibitem{RS03}
E.~J. Rawdon and R.~G. Scharein, \emph{Equilateral stick numbers}, to appear
  Contemp. Math., 2003.

\bibitem{SS02}
K.~P. Scannell and D.~P. Sinha, \emph{A one-dimensional embedding complex}, J.
  Pure Appl. Algebra \textbf{170} (2002), no.~1, 93--107.

\bibitem{Si02a}
D.~P. Sinha, \emph{Manifold-theoretic compactifications of configuration
  spaces}, In preparation.

\bibitem{Si00}
\bysame, \emph{The topology of spaces of knots}, Submitted.

\bibitem{St98}
T.~B. Stanford, \emph{Vassiliev invariants and knots modulo pure braid
  subgroups}, math.GT/9805092, 1998.

\bibitem{St70}
J.~D. Stasheff, \emph{${H}$-spaces from a homotopy point of view}, Lecture
  Notes in Math., vol. 161, Springer-Verlag, New York-Berlin-Heidelberg, 1970.

\bibitem{Va92}
V.~A. Vassiliev, \emph{Complements of discriminants of smooth maps: topology
  and applications}, Transl. Math. Monographs, vol.~98, Amer. Math. Soc.,
  Providence, 1992.

\bibitem{Vo03}
I.~Volic, \emph{Finite type knot invariants and calculus of functors}, Ph.D.
  thesis, Brown University, 2003.

\bibitem{We96}
M.~S. Weiss, \emph{Calculus of embeddings}, Bull. Amer. Math. Soc. (N.S.)
  \textbf{33} (1996), no.~2, 177--187.

\bibitem{We99}
\bysame, \emph{Embeddings from the point of view of immersion theory. {I}},
  Geom. Topol. \textbf{3} (1999), 67--101.

\end{thebibliography}
\end{document}